\documentclass{article}
\usepackage{amsmath, amssymb, amsthm}
\usepackage{enumitem}
\usepackage{bbm}  
\usepackage{thmtools, cleveref}
\usepackage{geometry}
\title{Equations driven by fast-oscillating functions of an It\^o diffusion process}
\author{
Tanner Reese
\footnote{Corresponding Author}
\footnote{Department of Mathematics, University of Arizona, Tucson, Arizona, United States}
\quad and \quad Jan Wehr
\footnote{Department of Mathematics and Program in Applied Mathematics,
University of Arizona, Tucson, Arizona, United States}
}

\newcommand\dr{\; \mathrm{d}}
\newcommand\eps[1]{{#1^{(\epsilon)}}}
\newcommand\domain[1]{\mathcal{D}\left(#1\right)}

\newcommand\probWhen[2]{\mathbb{P}\left(#1\,\middle|\,#2\right)}
\newcommand\expect[1]{\mathbb{E}\left[#1\right]}
\newcommand\expectWhen[2]{\mathbb{E}\left[#1\middle|#2\right]}

\newcommand\wbr[1]{\overline{#1}}

\newtheorem{defn}{Definition}

\newtheorem{lem}[defn]{Lemma}
\newtheorem{coro}[defn]{Corollary}
\newtheorem{prop}[defn]{Proposition}
\newtheorem{thm}[defn]{Theorem}

\begin{document}

\maketitle
\begin{abstract}
We study It\^o SDE systems driven by oscillating functions of a single It\^o diffusion process.
In the limit when oscillations become fast,
we show that the solution process converges in law to the process defined by an SDE
system driven by a multivariate Wiener process whose covariance we calculate explicitly.
Interestingly, the limiting system of SDEs are most naturally stated using the Stratonovich integral.
The problem has been originally motivated by experimental work and special cases of theorems proved here
provide a rigorous treatment of equations arising from physics.
\end{abstract}

\section{Introduction}
As a motivating example, for a real $ \epsilon > 0 $, a Wiener process $ W $, and
a $ 2\pi $-periodic function $ \phi : \mathbb{R} \to \mathbb{R} $ with $ \int_0^{2\pi} \phi \dr\theta = 0 $,
we define $ \eps{N}(t) := \frac{1}{\epsilon} \phi\left(\frac{1}{\epsilon} W(t)\right) $.
As $ \epsilon \downarrow 0 $, the law of $ \eps{N} $ approaches that of a white noise.
Abusing notation, we might say $ \eps{N}(t) \,dt \to c \dr \wbr{W}(t) $ in law as $ \epsilon \downarrow 0 $
where $ \wbr{W} $ is a Wiener process.
This fact can be deduced from the results in section 2.8.2 of \cite{KLO}, using the Wiener scaling.
It will also be derived as a very special case of the results presented here.
The present work studies stochastic differential equations (SDE) arising as limits of systems driven by such processes.
As far as the authors know, this has not been considered earlier.

More generally, suppose $ \{\phi_\alpha\}_{\alpha=1}^{d_B} $ is a collection of $ 2\pi $-periodic real functions
with antiderivatives $ \{\Phi_\alpha\}_{\alpha=1}^{d_B} $ so $ \Phi_\alpha' = \phi_\alpha $. 
Assume that the mean of $\phi_{\alpha}$ is zero: $ \int_0^{2\pi} \phi_{\alpha} \dr\theta = 0 $.
Then, $\Phi_{\alpha}$ is also periodic and can be chosen to have mean zero as well.  
If $ \eps{X} $ is a process in $ \mathbb{R}^{d_X} $ which satisfies the randomly driven equation
\begin{equation} \label{eq:wiener-sde}
\dr \eps{X}(t) = b\left( \eps{X}(t) \right) \dr t + \sum_{\alpha=1}^{d_B}
v_\alpha\left( \eps{X}(t) \right) \frac{1}{\epsilon} \phi_\alpha\left( \frac{1}{\epsilon} W(t) \right) \dr t
\end{equation}
then $ \eps{X} $ will converge in law to an It\^o diffusion $ X $.
It turns out that the equation describing the limiting process is most naturally written 
as a Stratonovich SDE:
\begin{equation} \label{eq:wiener-sde-limit}
\dr X(t) = b(\eps{X}(t)) \dr t + 2 \sum_{\alpha=1}^{d_B} v_{\alpha}(X(t)) \circ \dr B_\alpha(t)
\end{equation}
where $ \{B_\alpha\}_{\alpha=1}^{d_B} $ is a multivariate Wiener process with
\[
\expect{B_\alpha(s) B_\beta(t)} = \langle \Phi_\alpha, \Phi_\beta \rangle \min(s, t)
= \frac{\min(s, t)}{2\pi} \int_0^{2\pi} \Phi_\alpha(\theta) \Phi_\beta(\theta) \dr\theta.
\]
That is, the covariance matrix of $ (B_1, \ldots, B_d) $ is 
the Gram matrix of $ (\Phi_1, \ldots, \Phi_d) $ in $ L^2([0, 2\pi]) $ (with the normalized Lebesgue measure).
We stress that, while the original equation is driven by a function of a single Wiener process,
the limiting system is driven by a multivariate Wiener process of arbitrarily large dimension.
We refer to this behavior as the splitting of
a one-dimensional Wiener process into a multi-dimensional Wiener process.

Since, as stated above, the processes $ \frac{1}{\epsilon} \phi_\alpha\left( \frac{1}{\epsilon} W \right)$
approximate (and, therefore, can be thought of as a regularization of)
components of a white noise process, our results are related to the well-known Wong-Zakai theorem \cite{WZ}.
In this theorem, the Wiener process driving an SDE is approximated pathwise by a piecewise smooth process.
The solutions of equations with such regularized noise are shown to converge strongly to the solutions of the corresponding Stratonovich equations.
This suggests an alternate proof of (most of) our results.
First, one would prove convergence of the time integrals of the processes $ t \mapsto \frac{1}{\epsilon} \phi_\alpha\left( \frac{1}{\epsilon} W(t) \right) $
to components of a Wiener process.
Then, one would use Skorokhod's representation theorem to construct an almost surely convergent sequence
to which the Wong-Zakai theorem applies.
Instead, we find it more natural to prove the convergence in law of the solutions,
using convergence of the generators of the corresponding Markov semigroups.

In what follows, we replace the Wiener process $ W $ by a function of an It\^o diffusion.
Suppose $ M $ is a process in $ \mathbb{R}^{d_M} $
and $ W $ is a standard $ d_W $-dimensional Wiener process with independent components $ \{W_i\} $
such that
\begin{equation} \label{eq:driving-proc}
\dr M(t) = \nu(M(t)) \dr t + \sum_{i=1}^{d_W} \sigma_i(M(t)) \dr W_i(t)
\quad \text{for all} \;\; t \geq 0.
\end{equation}
We remark that the role played by the $ \{W_i\} $ is different here from that of a scalar Wiener process in \cref{eq:wiener-sde} -- the latter is replaced by a function of $ M $,
as we are going to describe in detail below.
In \cref{sect:results}, we use $ M $ to define two models studied in this paper: amplitude scaling and time scaling.
They are equivalent in the case when $ M $ is a Wiener process, but different in the general case.

\subsection{Motivation}
We briefly mention two mathematical models
which originally motivated this work.
Both came from experiments studying collective behavior of interacting agents.
One investigated aggregation and disaggregation in a system of communicating robots.
The other studied spontaneous formation of high-density regions in bacterial colonies.
Both came from the joint projects of one of the authors with experimental groups.
In each case, the theoretical part of the project relied on a detailed study
of one-particle (one-robot or one-bacterium) motion,
using it as a basis for a mean-field theory of the collective behavior.
A very particular case of the present work makes the one-particle calculation rigorous.

\subsubsection{Motion of a Phototactic Robot}
\label{ssub:robots}
The work \cite{MMWV} studies
a small robot moving in a plane,
randomly changing its direction and adapting its speed 
to its position, with a sensorial delay.
The equations describing the robot's motion are
\[ \begin{split}
\frac{\dr X_1}{\dr t} &= \frac{1}{\epsilon} u(X(t - k \epsilon^2)) \cos\left(\frac{1}{\epsilon} W(t)\right) \\
\frac{\dr X_2}{\dr t} &= \frac{1}{\epsilon} u(X(t - k \epsilon^2)) \sin\left(\frac{1}{\epsilon} W(t)\right)
\end{split} \]
where $ X = (X_1, X_2) $,
$ u $ is a position-dependent speed function,
$ W $ is a standard Wiener process,
and $ \epsilon $ is a scaling parameter.
One may attempt to analyze the limiting system as $ \epsilon \downarrow 0 $
by considering the linearization given by
\begin{equation} \begin{split} \label{eq:robot-linearized}
\dr X_1(t) &= \left[
-ku(X) \partial_1 u(X) \cos^2(\Theta)
- ku(X) \partial_2 u(X) \cos(\Theta) \sin(\Theta)
+ \frac{1}{\epsilon} u(X) \cos(\Theta)
\right] \dr t \\
\dr X_2(t) &= \left[-k u(X) \partial_1 u(X) \cos(\Theta) \sin(\Theta)
- ku(X) \partial_2 u(X) \sin^2(\Theta)
+ \frac{1}{\epsilon} u(X) \sin(\Theta)
\right] \dr t \\
\dr \Theta(t) &= \frac{1}{\epsilon} \dr W(t).
\end{split} \end{equation}
It is obtained by first linearizing $ X_1 $ and $ X_2 $ and
then linearizing $ u $ so as to get rid of the delay in the argument of the process $X$.
The resulting system of linear equations is then solved for $ \frac{\dr X_1}{\dr t} $ and $ \frac{\dr X_2}{\dr t} $.
A further approximation of matrix inversion leads to \cref{eq:robot-linearized}.
The robot is programmed to follow a discrete-time version of \cref{eq:robot-linearized}
with values of $ k $ which can be varied,
taking both positive and negative values.
\smallskip

In \cite{MMWV}, it is shown that in the limit $ \epsilon \downarrow 0 $,
the solution process of \cref{eq:robot-linearized} converges in law to the solution process of the system
\begin{equation} \label{eq:robot-limit} \begin{split}
\dr X_1(t) &= \left(1 - \frac{k}{2}\right) u(X) \partial_1 u(X) \dr t
+ \sqrt{2} u(X) \dr W_1(t) \\
\dr X_2(t) &= \left(1 - \frac{k}{2}\right) u(X) \partial_2 u(X) \dr t
+ \sqrt{2} u(X) \dr W_2(t)
\end{split} \end{equation}
where $ W_1 $ and $ W_2 $ are independent Wiener processes.
The noise terms in \cref{eq:robot-limit} come
from $ \frac{1}{\epsilon} \cos(\Theta) $ and $ \frac{1}{\epsilon} \sin(\Theta) $ in \cref{eq:robot-linearized}
since the two-dimensional process
\begin{equation} \label{eq:wiener-trigs}
\frac{1}{\sqrt{2} \epsilon} \left(
\int_0^t \cos\left(\frac{1}{\epsilon} W(s)\right) \dr s,
\int_0^t \sin\left(\frac{1}{\epsilon} W(s)\right) \dr s
\right)
\end{equation}
converges in law to the two-dimensional Wiener process.
This can be deduced from more general results of \cite{KLO} using Wiener scaling
or shown directly using the It\^o integration by parts formula
and the L\'evy characterization theorem.
Both the convergence in law of \cref{eq:wiener-trigs} to a two-dimensional Wiener process
and the convergence in law of the solutions of \cref{eq:robot-linearized}
to the solutions of \cref{eq:robot-limit}
follow from \cref{thm:ampl-scaling-limit}.

The origin of the terms $ - \frac{k}{2} u(X) \partial_1 u(X) $ and $ - \frac{k}{2} u(X) \partial_2 u(X) $
is also easy to explain
since, as $ \epsilon \downarrow 0 $,
the fast oscillations of the trigonometric functions of $ \frac{1}{\epsilon} W $
make $ \sin^2\left(\frac{1}{\epsilon} W\right) $ and $ \cos^2\left(\frac{1}{\epsilon} W\right) $
average to $ \frac{1}{2} $.
The appearance of the additional terms $ u(X) \partial_1 u(X) $
and $ u(X) \partial_2 u(X) $ can be explained by invoking the Wong-Zakai theorem;
since \cref{eq:wiener-trigs} approximates the two-dimensional Wiener process
by processes with smooth trajectories,
the drift of the limiting system must be corrected by a Stratonovich term.

\Cref{eq:robot-limit} which describes the motion of a single robot motion
in a time-independent speed landscape
is further used in \cite{MMWV} to study
the collective behavior of a system of several robots,
using the one-robot calculation as a mean-field approximation.
This leads to a prediction of a dynamical phase transition
from aggregating to disaggregating behavior
as the value of the parameter $ k $ in \cref{eq:robot-linearized} crosses $ -2 $.
This step is not further elaborated on in the present work.

\subsubsection{Motility Induced Phase Separation}
\label{ssub:mips}

In \cite{Stenh}, a model similar to \cref{eq:robot-linearized} is used to describe
the motion of a single bacterium in an environment created by other bacteria in a colony.
We focus here on the one-bacterium equations of motion.
These equations can be used with mean-field arguments to explain 
the spontaneous formation of regions of higher density.
This phenomenon is called MIPS (Motility Induced Phase Separation).

Motion of a single bacterium is described by a system of stochastic delay equations
which, in a linear approximation, can be replaced by the system
\begin{equation} \label{eq:mips} \begin{split}
\dr X_1(t) &= \left[\kappa \partial_1 w(X) \cos^2(\Theta)
+ \kappa \partial_2 w(X) \cos(\Theta) \sin(\Theta)
+ \frac{1}{\epsilon} w(X) \cos(\Theta) \right] \dr t \\
\dr X_2(t) &= \left[\kappa \partial_1 w(X) \cos(\Theta) \sin(\Theta)
+ \kappa \partial_2 w(X) \sin^2(\Theta)
+ \frac{1}{\epsilon} w(X) \sin(\Theta) \right] \dr t \\
\dr \Theta(t) &= \frac{1}{\epsilon} \dr W(t)
\end{split} \end{equation}
which belongs to the class of systems studied in the current paper.
By applying \cref{thm:ampl-scaling-limit},
one sees that the process defined by \cref{eq:mips},
converges in law to the solution of the It\^o SDE:
\[ \begin{split}
\dr X_1(t) &= \frac{1}{2} \kappa \partial_1 w(X)
+ w(X) \partial_1 w(X) + \sqrt{2} w(X) \dr W_1(t) \\
\dr X_2(t) &= \frac{1}{2} \kappa \partial_2 w(X)
+ w(X) \partial_2 w(X) + \sqrt{2} w(X) \dr W_2(t)
\end{split} \]
driven by two independent Wiener processes
or, equivalently, the Stratonovich SDE:
\[ \begin{split}
\dr X_1(t) &= \frac{1}{2} \kappa \partial_1 w(X)
+ \sqrt{2} w(X) \circ \dr W_1(t) \\
\dr X_2(t) &= \frac{1}{2} \kappa \partial_2 w(X)
+ \sqrt{2} w(X) \circ \dr W_2(t).
\end{split} \]

The derivation of the results in \cref{ssub:robots} and \cref{ssub:mips}
applies a multiscale ansatz to the Kolmogorov equation
and it is thus not entirely rigorous.
The present work provides its complete proof,
as a very special case of a much more general result.
The SDE system studied here can be thought of as
generalizing \cref{eq:robot-linearized} in several directions.
In the simplest case, we consider an SDE system not for two,
but for an arbitrary number of processes,
driven by a single one-dimensional Wiener process.
We replace the $ \sin $ and $ \cos $ functions
by arbitrary periodic functions of $ \theta $
and in the limit obtain a system driven by $ n $ independent Wiener processes.
We call this ``Wiener splitting".
Again, the drift terms are averaged over fast oscillations
and contain Stratonovich corrections.
We also consider systems driven by more general diffusion processes
which have more complicated limits,
as explained in detail in \cref{sect:results}.

\section{Results}
\label{sect:results}

\subsection{Amplitude Scaling}
\label{sub:ampl-scaling}

For a $ C^2 $-function $ \vartheta : \mathbb{R}^{d_M} \to \mathbb{R} $,
we replace $ W $ in \cref{eq:wiener-sde} by $ \vartheta(M) $
where $ M $ is as in \cref{eq:driving-proc}.
We allow the drift $ b $ to explicitly depend on the noise, obtaining
\begin{equation} \label{eq:ampl-scaling}
\dr \eps{X}(t) = b\left(\eps{X}(t), \frac{1}{\epsilon} \vartheta(M(t)) \right) \dr t
+ \frac{1}{\epsilon} \sum_{\alpha=1}^{d_B} v_\alpha\left(\eps{X}(t)\right)
\phi_\alpha\left(\frac{1}{\epsilon} \vartheta(M(t))\right) \dr t
\end{equation}
where $ b : \mathbb{R}^{d_X} \times \mathbb{R} \to \mathbb{R}^{d_X} $
and for all $ \alpha $, $ v_\alpha : \mathbb{R}^{d_X} \to \mathbb{R}^{d_X} $ and $ \phi_\alpha : \mathbb{R} \to \mathbb{R} $.
The main result of this section is that, under mild assumptions,
the solutions $ \eps{X} $ converge in law
to an It\^o diffusion satisfying an SDE similar to \cref{eq:wiener-sde-limit}.

\begin{thm} \label{thm:ampl-scaling-limit}
Suppose $ \{\phi_\alpha\}_{\alpha=1}^{d_B} \subseteq C_b(\mathbb{R}) $
and $ \{v_\alpha\}_{\alpha=1}^{d_B} \subseteq C_b^3(\mathbb{R}^{d_X}; \mathbb{R}^{d_X}) $
while $ b \in C_b^2(\mathbb{R}^{d_X} \times \mathbb{R}) $.
Additionally, $ \theta \mapsto b(x, \theta) $ and $ \{\phi_\alpha\} $ are $ 2\pi $-periodic
with $ \int_0^{2\pi} \phi_\alpha(\theta) \dr\theta = 0 $.
For every $ x \in \mathbb{R}^{d_X} $ and $ m \in \mathbb{R}^{d_M} $, we define
\[
\wbr{b}(x) := \frac{1}{2\pi} \int_0^{2\pi} b(x, \theta) \dr\theta
\quad \text{ and } \quad
\kappa(m) := \sqrt{\sum_{k=1}^{d_W} (\sigma_k(m) \cdot \nabla_m \vartheta(m))^2}
\]
and we assume that there exists $ \kappa_0 > 0 $
such that for all $ t > 0 $, $ \kappa(M(t)) > \kappa_0 $ almost surely.
It follows from the \cref{lem:periodic-diff}  below that for each $ \alpha $, there exists a mean-zero continuous $ 2\pi $-periodic function $ \Phi_\alpha $ with $ \Phi_\alpha' = \phi_\alpha $.
Then, the solution process $ \eps{X}(t) $ of the amplitude scaling equation (\ref{eq:ampl-scaling})
converges in law to the solution of the SDE
\begin{equation} \label{eq:ampl-limit-sde} \begin{split}
\dr X(t) &= \wbr{b}(X(t)) \dr t + \frac{2}{\kappa(M(t))} \sum_{\alpha=1}^{d_B} v_\alpha(X(t)) \circ \dr \tilde{B}_\alpha(t) \\
\dr M(t) &= \nu(M(t)) \dr t + \sum_{k=1}^{d_W} \sigma_k(M(t)) \dr W_k(t)
\end{split} \end{equation}
where $ \{\tilde{B}_\alpha\}_{\alpha=1}^{d_B} $ is a multivariate Wiener process
independent of $ \{W_k\}_{k=1}^{d_W} $
with $ \expect{\tilde{B}_\alpha(s) \tilde{B}_\beta(t)} = \frac{1}{2\pi} \langle \Phi_\alpha, \Phi_\beta \rangle_{L^2([0, 2\pi))} \min(s, t) $.
\end{thm}

\subsection{Time Scaling}

Alternatively, we may replace $ \phi\left( \frac{1}{\epsilon} W(t) \right) $ (which has the same law as $ \phi\left(\big. W(t / \epsilon^2) \right) $) by $ \phi(M(t / \epsilon^2)) $ to get
\begin{equation} \label{eq:time-scaling}
\dr \eps{X}(t) = b\left(\eps{X}(t), M(t / \epsilon^2)\right) \dr t + \frac{1}{\epsilon} \sum_{\alpha=1}^{d_B} 
v_\alpha\left(\eps{X}(t)\right) \phi_\alpha(M(t / \epsilon^2)) \dr t
\end{equation}
which we will call the time scaling.
The limiting behavior of this system is more difficult to establish than that of the amplitude scaling above.
In particular, it requires additional assumptions on the process $ M $.

\subsubsection{Time Scaling in a Non-Compact Space using Ergodicity}
\label{ssub:ergodic-time-scale}

One option to ensure that $ \eps{X} $ converges in law as $ \epsilon \downarrow 0 $
is to require $ M $ to be strongly ergodic
so that it converges sufficiently quickly to its invariant measure.

\begin{defn}
We say that a stochastic process $ M $ with values in a measurable space $ (S, \mathcal{M}) $
is \textbf{strongly ergodic} if
there exist a probability measure $ \mu $ on $ S $
and real numbers $ C, \lambda, t_0 > 0 $ such that
\[
\sup_{E \in \mathcal{M}} \left|
\probWhen{M(t) \in E \Big.}{M(0) = m} - \mu(E)
\right| \leq C e^{-\lambda t}
\quad \text{for all} \quad
m \in S, \; E \in \mathcal{M}, \text{ and } t \geq t_0.
\]
\end{defn}

The strong ergodicity of $ M $ will ensure the existence of a restricted inverse for the infinitesimal generator.
This inverse will be crucial to prove the convergence in law of \cref{eq:time-scaling}.
We will consider the case where $ S $ is a locally compact topological space
with Borel $ \sigma $-algebra $ \mathcal{M} $
and $ M $ is a Feller-Dynkin process (in the sense of \cite{RW}) on $ S $.
For any $ f \in C_b(S) $ (i.e. bounded continuous functions on $ S $), we have $ T_M(t) f \in C_b(S) $
given by $ (T_m(t) f)(m) := \expectWhen{f(M_t)}{M_0 = m} $ for all $ m \in S $.
Following \cite{Dynk}, we define $ \mathcal{B}_M $ as
\[ \mathcal{B}_M := \left\{f \in C_b(S) \,:\, \lim_{t \downarrow 0} T_M(t) f = f \right\} \]
which is a Banach space by 1.3B of \cite{Dynk}.
Also, by 1.3A of \cite{Dynk}, the semigroup $ T_M $ is strongly continuous on $ \mathcal{B}_M $.

\begin{lem} \label{lem:ergodic-inv}
Suppose $ M $ is a strongly ergodic Feller-Dynkin process on $ (S, \mathcal{M}) $
with constants $ C, \lambda, t_0 > 0 $.
If $ A_M $ is the infinitesimal generator for $ T_M $
with $ \domain{A_M} \subseteq \mathcal{B}_M $
and $ \mu $ is a probability measure invariant under the action of the semigroup of $ M $
then there exists $ A_M^{-1} : D_0 \to D_0 $
where $ D_0 := \{ \phi \in \mathcal{B}_M \,:\, \int \phi \dr\mu = 0\} $
such that $ A_M A_M^{-1} = I_{D_0} $
and $ \|A_M^{-1}\| \leq t_0 + \frac{C}{\lambda} $.
\end{lem}
In order to state our result in a form analogous to \cref{thm:ampl-scaling-limit},
we introduce the following bilinear form on the space $ D_0 $.

\begin{defn} \label{defn:inv-gener-form}
Suppose $ M $ is a strongly ergodic Feller-Dynkin process with invariant probability measure $ \mu $.
We define the bilinear form associated with $ M $ on $ D_0 $ as
\[
\langle \phi, \psi \rangle_M := -\int_S \phi(m) \; (A_M^{-1} \psi)(m) \dr\mu(m)
= -\int_0^\infty \int_S \phi \cdot T_M(t) \psi \dr\mu \dr t
\]
for all $ \phi, \psi \in D_0 $
where $ A_M^{-1} $ is the same as in \cref{lem:ergodic-inv} and $T_M$ is the semigroup for $ M $.
\end{defn}

We will also define $ \wbr{b} : \mathbb{R}^{d_X} \to \mathbb{R}^{d_X} $ as
\[ \wbr{b}(x) := \int_S b(x, m) \dr \mu(m) \]
where we assume that for every $ x \in \mathbb{R}^{d_X} $, $ m \mapsto b(x, m) $ is $ \mu $-integrable.

\begin{thm} \label{thm:ergodic-time-scaling}
Suppose $ S $ is a locally compact topological space,
$ \mathcal{M} $ is its Borel $ \sigma $-algebra,
and $ M $ is a strongly ergodic Feller-Dynkin process on $ (S, \mathcal{M}) $.
We assume $ \mathcal{B}_M $ is an algebra,
$ \{\phi_\alpha\}_{\alpha=1}^{d_B} \subseteq \mathcal{B}_M $ are real,
and $ \{v_\alpha\}_{\alpha=1}^{d_B} \subseteq C_b^3(\mathbb{R}^{d_X}; \mathbb{R}^{d_X}) $
while $ b \in C_b(R; \mathcal{B}_M) $ is second-differentiable in $ x $.
For every $ \epsilon > 0 $, let $ \eps{X} $ be the solution of the time scaling \cref{eq:time-scaling}.
We assume that the restriction of $ \langle \cdot, \cdot \rangle_M $
to the subspace $ \mathrm{Span}\{\phi_\alpha\}_{\alpha=1}^{d_B} $
is symmetric and positive semi-definite.
Then, the solutions $ \eps{X} $ converge in law to the solution of the SDE
\begin{equation} \label{eq:time-limit-sde} \begin{split}
\dr X(t) = \wbr{b}(X(t)) \dr t + \sqrt{2} \sum_{\alpha=1}^{d_B} v_\alpha(X(t)) \circ \dr \tilde{B}_\alpha(t)
\end{split} \end{equation}
where $ \{\tilde{B}_\alpha\}_{\alpha=1}^{d_B} $ is a multivariate Wiener process
with $ \expect{\tilde{B}_\alpha(s) \tilde{B}_\beta(t)} = \langle \phi_\alpha, \phi_\beta \rangle_M \min(s, t) $
for $ s, t \geq 0 $.
\end{thm}

As shown in Wang \cite{JW},
any process confined by a sufficiently strong potential
will be strongly ergodic.
Thus, we can apply \cref{thm:ergodic-time-scaling} to a broad class of processes. For a twice continuously differentiable potential $ U : \mathbb{R}^{d_M} \to \mathbb{R} $ we define $ V := \frac{1}{2} (|\nabla U|^2 - \nabla^2 U) $.
We require that
\begin{equation} \label{eq:poten-conds} \begin{split}
\text{there exist } a, c > 0 &\text{ such that for all } z, \; U(z) \geq a |z| - c \\
\text{and there exists } c_V > 0 &\text{ such that for all } z, \;
V(x) > -c_V \text{ and } V(z) \to \infty \text{ as } |z| \to \infty.
\end{split} \end{equation}

\begin{coro} \label{coro:ergodic-wiener-well}
Suppose $ M $ is a solution to the SDE
\[ \dr M(t) = - (\nabla_m U)(M(t)) \dr t + \dr W(t) \]
where $ U : \mathbb{R}^{d_M} \to \mathbb{R} $ fulfills \cref{eq:poten-conds}
and $ W $ is a standard Wiener process in $ \mathbb{R}^{d_M} $.
We assume there exist differentiable $ h_1, h_2 : \mathbb{R}_+ \to \mathbb{R}_+ $ with $ h_1(0) = h_2(0) = 0 $ such that
\begin{equation} \label{eq:wang-conds} \begin{split}
h_1'(|x|) \leq \xi \cdot (\nabla \nabla U) \xi \leq h_2'(|x|)
\quad &\text{ for all } x, \xi \in \mathbb{R}^{d_M}
\text{ with } |\xi| = 1, \\
\int_{r_0}^\infty \frac{\dr r}{r h_1(r)} < \infty
\quad \text{ for some } r_0 > 0,
&\text{ and } \quad \lim_{r \to \infty} h_1(r) \int_r^\infty \frac{\dr s}{s h_2(s)} = \infty
\end{split} \end{equation}
where $ \nabla \nabla U $ is the Hessian matrix of $ U $.
Additionally, we assume $ b \in C_b^2(\mathbb{R}^{d_X} \times \mathbb{R}^{d_M}) $
with $ (\nabla_m U) \cdot (\nabla_m b) \in C_b(\mathbb{R}^{d_X} \times \mathbb{R}^{d_M}) $
and for all $ \alpha $, $ v_\alpha \in C_b^3(\mathbb{R}^{d_X}; \mathbb{R}^{d_X}) $
and $ \phi_\alpha \in C_b^2(\mathbb{R}^{d_M}) $
with $ (\nabla_m U) \cdot (\nabla_m \phi_\alpha) \in C_b(\mathbb{R}^{d_M})$.
If $ \eps{X} $ are solutions to \cref{eq:time-scaling}
then $ \eps{X} \to X $ in law where $ X $ is a solution of
\[
\dr X(t) = \wbr{b}(X(t)) \dr t + \sqrt{2} \sum_{\alpha=1}^{d_B} v_\alpha(X(t)) \circ \dr \tilde{B}_\alpha(t)
\]
in which, using the operator $ H := -\frac{1}{2} \nabla_m^2 + V = -\frac{1}{2} \left(\nabla_m^2 + \nabla_m^2 U - |\nabla_m U|^2 \right) $,
the covariance of $ \{\tilde{B}_\alpha\}_{\alpha=1}^{d_B} $ is
\[
\expect{\tilde{B}_\alpha(s) \tilde{B}_\beta(t)} =
\min(s, t) \frac{1}{\int e^{-2U} \dr m}\int (e^{-U} \phi_\alpha)(m) \cdot (H^{-1} e^{-U} \phi_\beta)(m) \dr m.
\]
\end{coro}

Without \cref{eq:wang-conds}, even if $ U $ grows sufficiently quickly,
if its growth is not consistently quick,
$ M $ may be prevented from filling certain areas of $ \mathbb{R}^{d_M} $ fast enough
precluding strong ergodicity.

\subsubsection{Time Scaling driven by Integrated Noise}

In applications, the driving process may be defined
as a time integral of another (typically stationary) process. 
A standard example:  it is well known that if $ Z $ is an Ornstein-Uhlenbeck (OU) process
(i.e. a diffusion in the well $ U(z) = \frac{1}{2} C z^2 $, see \cref{eq:integ-noise-driving})
then $ t \mapsto \frac{1}{\epsilon} Z(t / \epsilon^2) $ 
approximates a white noise as $ \epsilon \downarrow 0 $ (theorem 9.5 of \cite{Nels}).
Hence, the integral $ t \mapsto \frac{1}{\epsilon}\int_0^t Z(s / \epsilon^2) \dr s $
approximates a Wiener process.
So, if $ \Theta(t) := \int_0^t Z(s) \dr s $
then $\Theta(t / \epsilon^2) = \frac{1}{\epsilon^2} \int_0^t Z(s / \epsilon^2) \dr s $
may replace $ M(t / \epsilon^2) $ as the driving process in  \cref{eq:time-scaling}.
We will generalize this case
by taking $ Z $ to be a diffusion process in a potential belonging to a wide class 
and $ \Theta $ to be the integral of some function of $ Z $:
\begin{equation} \label{eq:integ-noise-driving}
\dr \Theta(t) = \rho(Z(t)) \dr t
\quad \text{and} \quad
\dr Z(t) = -\nabla U(Z(t)) \dr t + \dr W(t)
\end{equation}
where $ \rho, U : \mathbb{R}^{d_Z} \to \mathbb{R} $.
Then, after scaling time by $ \epsilon $,
we obtain the following time-scaling SDE for $ \left(\eps{X}, \eps{\Theta}, \eps{Z}\right) $
\begin{equation} \label{eq:integ-noise-time-scaling} \begin{split}
\dr \eps{X}(t) &= b\left(\eps{X}(t), \eps{\Theta}(t)\right) \dr t
+ \frac{1}{\epsilon} \sum_{\alpha=1}^{d_B} v_\alpha\left(\eps{X}(t)\right) \phi_\alpha\left(\eps{\Theta}(t)\right) \dr t \\
\dr \eps{\Theta}(t) &= \frac{1}{\epsilon^2} \rho\left( \eps{Z}(t) \right) \dr t \\
\dr \eps{Z}(t) &= -\frac{1}{\epsilon^2} (\nabla U)\left( \eps{Z}(t) \right) \dr t + \frac{1}{\epsilon} \dr W(t).
\end{split} \end{equation}
As with the previous theorems,
there is a bilinear form $ \langle \cdot, \cdot \rangle_{\Theta,Z} $ defined in \cref{eq:integ-noise-form}
which describes the covariance matrix of the Wiener noise in the limit.

\begin{thm} \label{thm:integ-noise-time-scaling}
For every $ \epsilon > 0 $, suppose that $ \left( \eps{X}, \eps{\Theta}, \eps{Z} \right) $
is a solution of \cref{eq:integ-noise-time-scaling}.
We assume $ U $ fulfills \cref{eq:poten-conds}
and, additionally, $ |\rho|^2 \leq a_V (V + c_V) $ for some reals $ a_V, c_V > 0 $.
We also assume that for every $ \alpha $, $ v_\alpha \in C_b^3(\mathbb{R}^{d_X}; \mathbb{R}^{d_X}) $,
and $ \phi_\alpha $ is a trigonometric polynomial
and $ b $ can be written as
\begin{equation} \label{eq:integ-noise-drift-form}
b(x, m) := b_0(x) + \sum_{j=1}^k b_j(x) \psi_j(m)
\end{equation}
where $ \{\psi_j\}_{j=1}^k $ are trigonometric polynomials
and $ \{b_j\}_{j=1}^k \in C_b^2(\mathbb{R}^{d_X}; \mathbb{R}^{d_X}) $.
Then, the form $ \langle \cdot, \cdot \rangle_{\Theta,Z} $
given by \cref{eq:integ-noise-form}
for the driving process $ (\Theta, Z) $ (without time scaling)
is well-defined for trigonometric functions
and $ \eps{X} \to X $ in law as $ \epsilon \downarrow 0 $
where $ X $ is a solution to the SDE
\begin{equation} \label{eq:integ-noise-limit-sde}
\dr X(t) = b_0(X(t)) \dr t
+ \sqrt{2} \sum_{\alpha=1}^{d_B} v_\alpha(X(t)) \circ \dr \tilde{B}_\alpha(t)
\end{equation}
where $ \{\tilde{B}_\alpha\}_{\alpha=1}^{d_B} $ is a multivariate Wiener process
with $ \expect{\tilde{B}_\alpha(s) \tilde{B}_\beta(t)} = \langle \phi_\alpha, \phi_\beta \rangle_{\Theta,Z} \min(s, t) $ for $ s, t \geq 0 $.
\end{thm}

This theorem can be restricted to a case of particular interest.
As mentioned earlier, for $ U(z) = \frac{1}{2} z^2 $, $ Z $ becomes an Ornstein-Uhlenbeck (OU) process.
We can also take $ \rho(z) = z $
so that $ \Theta $ is simply the integral of $ Z $.
In this case, we can explicitly calculate
the values of the form $ \langle \cdot, \cdot \rangle_{\Theta,Z} $.

\begin{prop} \label{prop:ou-proc-form}
Suppose that the (unscaled) driving process \cref{eq:integ-noise-driving} as used in \cref{eq:integ-noise-time-scaling}
is the solution of the SDE
\[
\dr \Theta(t) = Z(t) \dr t
\quad \text{and} \quad
\dr Z(t) = - Z(t) \dr t + \dr W(t).
\]
Then, for any $ k, \ell \in \mathbb{Z} $,
the form $ \langle \cdot, \cdot \rangle_{\Theta,Z} $ from \cref{defn:inv-gener-form} is
\[
\langle e^{i \ell \theta}, e^{ik \theta} \rangle_{\Theta,Z}
= e^{\frac{k^2}{2}} \delta_{k\ell} \sum_{n=0}^\infty \frac{(-1)^n k^{2n}}{(n + k^2/2) 2^n n!}
\]
which for $ k = l =1 $ becomes $ \langle e^{i \theta}, e^{i \theta} \rangle_{\Theta,Z} = \sqrt{2e\pi} \, \mathrm{erf}\left(\frac{1}{\sqrt{2}}\right) $.
\end{prop}

\filbreak
\section{Proofs}

\subsection{Scaling Limit for General Processes}

Before we state and prove \cref{prop:semigrp-average}, we make several introductory comments. 
All of the cases dealt with in this text have a similar structure.
Suppose $ R $ and $ S $ are locally compact topological spaces.
The space $ S $ will be averaged out after taking the limit
while $ R $ will remain.
We restrict our attention to the case
where $ R = \mathbb{R}^d $ is a Euclidean space.
In the particular case of \cref{eq:wiener-sde},
the space $ S = \mathbf{S}^1 \cong \mathbb{R} / 2\pi \mathbb{Z} $
and for $ R $, $ d = d_X $ so $ R = \mathbb{R}^{d_X} $.

For each $ \epsilon > 0 $, let $ \{\eps{Y}(t)\}_{t \geq 0} $ be a Markov process on $ R \times S $.
We decompose it as $ \left(\eps{Q}, \eps{\Theta}\right) := \eps{Y} $
where $ \eps{Q} $ is a process on $ R $
and $ \eps{\Theta} $ is a process on $ S $.
Additionally, let $ \{Q(t)\}_{t \geq 0} $ be a Markov process on $ R $.
Our goal is to show that $ \eps{Q} \to Q $ in law
whenever the generator of $ \eps{Y} $ fulfills \cref{eq:eps-gener-defn}
and the generator of $ Q $ fulfills \cref{eq:limit-gener-defn}.
\big(For \cref{eq:wiener-sde},
we have $ \eps{Q} = \eps{X} $, $ Q = X $, and $ \eps{\Theta} = \frac{1}{\epsilon} W $.
However, in other cases, $ \eps{Q} $ may also include the driving process
and $ \eps{\Theta} $ becomes more complicated.\big)

We will argue on the basis of semigroups on an appropriate function space.
Let $ \mu $ be a Borel probability measure on $ S $
and $ \mathcal{B} $ be a Banach space of
(not necesarily all) $ \mu $-integrable complex functions on $ S $.
\big(For \cref{eq:wiener-sde}, $ \mu $ is the uniform probability measure on $ \mathbf{S}^1 $
and $ \mathcal{B} = C(\mathbf{S}^1) $.\big)
Now, let $ C_0(R; \mathcal{B}) $ be the space of continuous $ \mathcal{B} $-valued functions on $ R $
vanishing at infinity endowed with the uniform metric.
Then, suppose $ \left\{\eps{T_t}\right\}_{t \geq 0} $ is the strongly continuous semigroup
for $ \eps{Y} $ on $ C_0(R; \mathcal{B}) $
with infinitesimal generator $ \eps{A} $.
Also, suppose $ \{T_t\}_{t \geq 0} $ is the strongly continuous semigroup for $ Q $ on $ C_0(R) $
with infinitesimal generator $ A $.
We assume that for sufficiently regular $ f \in C_0(R; \mathcal{B}) $,
$ \eps{A} $ can be written as
\begin{equation} \label{eq:eps-gener-defn}
\eps{A} f = L f
+ \frac{1}{\epsilon} (E + J) f
+ \frac{1}{\epsilon^2} \eta K f
\quad \text{where} \quad
J f := \sum_{\alpha=1}^n \phi_\alpha v_\alpha \cdot \nabla_R f
= \sum_{\alpha=1}^n \phi_\alpha \nabla_{v_\alpha} f
\end{equation}
while $ E $ and $ L $ are densely-defined operators on $ C_0(R; \mathcal{B}) $,
$ \eta \in C_b(R) $ is invertible in $ C_b(R) $,
$ K $ is a densely-defined operator on $ \mathcal{B} $,
$ \{\phi_\alpha\}_{\alpha=1}^n \subseteq \mathcal{B} $,
and $ \{v_\alpha\}_{\alpha=1}^n \subseteq C_b(R; \mathbb{R}^d) $
are bounded $ \mathbb{R}^d $-valued functions.
Here, $ \nabla_R $ is the gradient on $ R = \mathbb{R}^d $
and $ \nabla_{v_\alpha} = v_\alpha \cdot \nabla_R $ is the directional derivative
with respect to the vector field $ v_\alpha $.

One may think of $ K $ as describing the ``fast" dynamics on $ S $
(of order $ 1 / \epsilon^2 $)
which may be affected by $ R $ via $ \eta $.
These dynamics will be averaged out in the limit.
\big(For \cref{eq:wiener-sde}, $ K = \frac{1}{2} \nabla_\theta^2 $ and $ \eta = 1 $.\big)
This averaging is described by
the continuous linear functional $ \omega : \mathcal{B} \to \mathbb{C} $ given by
\[
\omega(f) := \int_S f(s) \dr\mu(s)
\quad \text{for all} \quad f \in \mathcal{B}.
\]
In the limiting equation,
this averaging causes each $ \phi_\alpha\left(\eps{\Theta}\right) $
to become a white noise.
A critical step in our argument
is that appropriate functions in $ \ker(\omega) $
are in the image of $ K $.
Thus, ``variations of $ \mu $-mean zero" in $ S $
can be ``canceled out" by $ K $.
More precisely, for some linear subspace $ \Omega \subseteq \mathcal{B} $ with $ 1 \in \Omega $,
we will need $ K $ to be invertible on $ \ker(\omega) \cap \Omega $.
One may think of $ \Omega $ as a space of functions
on which $ K $ behaves ``nicely".
For most of the cases we consider,
$ \Omega $ will simply be $ \mathcal{B} $.
However, for \cref{thm:integ-noise-time-scaling} and \cref{prop:ou-proc-form},
$ \Omega $ will be restricted to trigonometric polynomials.
\smallskip

The operator $ L $ describes dynamics on $ R $
that will be averaged by $ \omega $ in the limit.
\big(For \cref{eq:wiener-sde}, $ L = b(x) \cdot \nabla_x $.\big)
The operator $ E $ describes the ``slow" dynamics on $ S $
(of order $ 1 / \epsilon $)
which may be $ R $-dependent.
The effect of $ E $ will disappear in the limit.
For most of the cases we consider,
excluding those in \cref{sub:ampl-scaling},
the operator $ E $ will be zero.
Intuitively, $ E $ has zero ``average" effect on the dynamics of $ S $.
\smallskip

For any $ \phi, \psi \in \mathcal{B} $,
if $ \psi \in \ker(\omega) \cap \Omega $
then, since $ K $ is invertible on $ \ker(\omega) \cap \Omega $, there exists $ K^{-1} \psi \in \mathcal{B} $
and, further, if $ \phi \cdot K^{-1} \psi \in \Omega $
then we may define the bilinear form
\begin{equation} \label{eq:gener-biform}
\langle \phi, \psi \rangle_K := -\omega(\phi^* \cdot K^{-1} \psi)
\end{equation}
where $ \phi^* $ is the complex conjugate.
We will require that for all $ \alpha $ and $ \beta $,
one has $ \phi_\alpha, \phi_\beta \in \ker(\omega) \cap \Omega $
and $ \phi_\alpha \cdot K^{-1} \phi_\beta \in \Omega $.
Then, the Gram matrix of $ \{\phi_\alpha\} $
with respect to $ \langle \cdot, \cdot \rangle_K $
will describe the covariance of the white noises in the limit.
We assume that for sufficiently regular $ f \in C_0(R) $, $ A $ can be written as
\begin{equation} \label{eq:limit-gener-defn}
A f = \omega(L f) + \eta^{-1} \sum_{\alpha,\beta=1}^n \langle \phi_\alpha, \phi_\beta \rangle_K \nabla_{v_\alpha} \nabla_{v_\beta} f.
\end{equation}
If $ \langle \cdot, \cdot \rangle_K $ is symmetric and positive semi-definite on $ \{\phi_\alpha\} $
then the matrix $ C = \{c_{\alpha \beta}\}_{\alpha,\beta=1}^n := \{\langle \phi_\alpha, \phi_\beta \rangle_K\}_{\alpha,\beta=1}^n $
is also symmetric and positive semi-definite.
Thus, there exists $ H = \{h_{\alpha \gamma}\}_{\alpha,\gamma=1}^n := \sqrt{C} $
and we can define vector fields $ \{\tilde{v}_\gamma\}_{\gamma=1}^n \subseteq C_b(R) $
as $ \tilde{v}_\gamma := \sum_{\alpha=1}^n h_{\alpha \gamma} v_\alpha $.
If $ \nabla_{v_\alpha} (\eta^{-1/2}) = 0 $ for all $ \alpha $ then
\begin{equation} \label{eq:decoupled-vec-fields} \begin{split}
Af &= \omega(L f) + \eta^{-1} \sum_{\gamma=1}^n \nabla_{\tilde{v}_\gamma} \nabla_{\tilde{v}_\gamma} f
= \omega(L f) + \eta^{-1} \sum_{\gamma=1}^n \left[\big.
(\nabla_{\tilde{v}_\gamma} \tilde{v}_\gamma) \cdot \nabla_R f
+ (\nabla \nabla_R f)(\tilde{v}_\gamma, \tilde{v}_\gamma)\right] \\
&= \omega(L f) + \sum_{\gamma=1}^n \left[
\frac{1}{2} \left(\nabla_{\sqrt{2 \eta^{-1}} \tilde{v}_\gamma} \sqrt{2 \eta^{-1}} \tilde{v}_\gamma\right) \cdot \nabla_R f
+ \frac{1}{2} (\nabla \nabla_R f)\left(\sqrt{2 \eta^{-1}} \tilde{v}_\gamma, \sqrt{2 \eta^{-1}} \tilde{v}_\gamma\right) \right].
\end{split} \end{equation}
where $ \nabla \nabla_R f $ is the Hessian of $ f $
treated as a quadratic form on $ \mathbb{R}^d $.
We define $ D = C_0^2(R) \subseteq C_0(R) $ as
\begin{equation} \label{eq:smooth-vanish-defn}
D := \{f \in C_0(R) = C_0(\mathbb{R}^d) \,:\, f \text{ is second-differentiable with its first and second derivatives in } C_0(R)\}.
\end{equation}

Suppose $ B $ is a Banach space
and $ F : \domain{F} \to B $ is a densely-defined operator
with $ \domain{F} \subseteq B $.
For any linear subspace $ V \subseteq B $,
we say $ F $ \textbf{preserves} $ V $
if $ F(V \cap \domain{F}) \subseteq V $.
For any linear subspace $ V \subseteq \domain{F} $,
we say $ F $ is \textbf{right-invertible on $ V $}
if there exists a linear operator $ F^{-1} : V \to V $
such that for all $ v \in V $, $ F F^{-1} v = v $. 
The above comments set the stage for the detailed statement and proof of the main technical result.
\medskip

\begin{prop} \label{prop:semigrp-average}
Suppose $ \mu $ is a Borel probability measure on $ S $,
$ \mathcal{B} $ is a Banach space of (not necessarily all) $ \mu $-integrable complex functions on $ S $,
and $ \Omega \subseteq \mathcal{B} $ with $ 1 \in \Omega $.
We write $ \omega(f) := \int f \dr\mu $ for all $ f \in \mathcal{B} $.
Also, suppose $ \eta \in C_b(R) $ is positive and invertible in $ C_b(R) $,
$ \{v_\alpha\}_{\alpha=1}^n \subseteq C_b^1(R; \mathbb{R}^d) $,
and $ \{\phi_\alpha\}_{\alpha=1}^n \subseteq \ker(\omega) \cap \Omega $ are real
while $ L $ and $ E $ are densely-defined operators on $ C_0(R; \mathcal{B}) $
and $ K $ is a densely-defined operator on $ \mathcal{B} $
such that
\vspace{-0.5em}
\begin{enumerate}[nosep]
\item $ K(1) = 0 $ and $ K $ is right-invertible on $ \ker(\omega) \cap \Omega $
\item $ \nabla_{v_\alpha} (\eta^{-1/2}) = 0 $
and $ \phi_\alpha \cdot K^{-1} \phi_\beta \in \Omega $ for all $ \alpha $ and $ \beta $
\item $ E $ is $ C_0(R) $-linear,
$ C_0(R) \subseteq \ker(E) $,
and $ E $ preserves $ \ker(\omega) \cap \Omega $.
\item $ D \subseteq \domain{L} $,
for all $ \phi \in \mathcal{B} $ and $ f \in \domain{L} $, $ L(\phi f) = \phi L(f) $,
and either $ L(D) \subseteq C_0(R) \otimes \Omega $
(where $ \otimes $ is the algebraic tensor product)
or $ K^{-1} $ is bounded on $ \ker(\omega) \cap \Omega $ and $ L(D) \subseteq C_0(R; \Omega) $.
\end{enumerate}
Let $ \eps{Y} = \left(\eps{Q}, \eps{\Theta}\right) $ and $ Q $ be Markov processes on $ R \times S $ and $ R $, respectively, as defined above
with infinitesimal generators $ \eps{A} $ and $ A $.
If $ \eps{A} $ fulfill \cref{eq:eps-gener-defn},
$ A $ fulfills \cref{eq:limit-gener-defn},
$ D $ is a core for $ A $,
and $ \langle \cdot, \cdot \rangle_K $
(as defined in \cref{eq:gener-biform})
is symmetric and positive semi-definite on $ \{\phi_\alpha\} $
then $ \eps{Q} \to Q $ in law.
\end{prop}
\begin{proof}
The convergence of $ \eps{Q} $ to $Q$ will follow from the convergence of the associated semigroups.
Because of theorem 6.1 in the first chapter of \cite{EK},
since $ D $ is a core for $ A $,
in order to prove uniform convergence of the semigroups on bounded intervals,
it suffices to show that for any $ f \in D $,
there exist $ \eps{f} \in \domain{\eps{A}} $ for all $ \epsilon > 0 $
such that $ \eps{f} \to f $
and $ \eps{A} \eps{f} \to A f $ as $ \epsilon \downarrow 0 $.
Given $ f \in D $, we will construct $ \eps{f} \in \domain{\eps{A}} $
as $ \eps{f} := f + \epsilon f_1 + \epsilon^2 f_2 $
with $ f_1, f_2 \in C_0(R; \mathcal{B}) $.
\smallskip

Because $ f \in D $,
for every $ \alpha $, we have $ \nabla_{v_\alpha} f \in C_0(R) $.
Also, $ K $ is right-invertible on $ \ker(\omega) \cap \Omega $
and $ \phi_\alpha \in \ker(\omega) \cap \Omega $
so $ K^{-1} \phi_\alpha \in \ker(\omega) \cap \Omega \subseteq \Omega $.
Then, $ \eta^{-1} \nabla_{v_\alpha} f \in C_0(R) $
so $ \eta^{-1} (K^{-1} \phi_\alpha) \nabla_{v_\alpha} f \in C_0(R) \otimes \Omega $
and we define
\begin{align}
\label{eq:first-ord-defn}
&f_1 := -\eta^{-1} K^{-1} J f = - \eta^{-1} \sum_{\alpha=1}^n (K^{-1} \phi_\alpha) 
\nabla_{v_\alpha} f \in C_0(R) \otimes \Omega \\
\label{eq:first-ord-vanish}
\text{and, since $ f \in C_0(R) \subseteq \ker(E) $,} \quad
&\eta K f_1 + (E + J) f = \eta K (- \eta^{-1} K^{-1} J f) + J f = - J f + J f = 0.
\end{align}
For all $ \alpha $, since $ K^{-1} \phi_\alpha \in \ker(\omega) \cap \Omega $
and $ E $ preserves $ \ker(\omega) \cap \Omega $,
we know $ E K^{-1} \phi_\alpha \in \ker(\omega) \cap \Omega $.
Also, for all $ \alpha $ and $ \beta $,
by assumption, $ \phi_\alpha \cdot K^{-1} \phi_\beta \in \Omega $
The second derivatives of $ f $
are continuous and vanish at infinity
and $ \{v_\alpha\} $ have bounded continuous derivatives
so $ f \in \domain{\nabla_{v_\alpha} \nabla_{v_\beta}} $.
Note that since $ \nabla_{v_\alpha} $ is a derivation, $ \nabla_{v_\alpha} (\eta^{-1/2}) = 0 $ 
implies that $ \nabla_{v_\alpha} $ commutes with multiplication by $ \eta^{-1/2} $.
Firt, assume that $ L f \in L(D) \subseteq C_0(R) \otimes \Omega $,
since $ E $ is $ C_0(R) $-linear,
we can define
\[
g := L f + (E + J) f_1
= L f - \eta^{-1} \sum_{\alpha=1}^n (E K^{-1} \phi_\alpha) \nabla_{v_\alpha} f
- \eta^{-1} \sum_{\alpha,\beta=1}^n (\phi_\alpha \cdot K^{-1} \phi_\beta)
\;\; \nabla_{v_\alpha} \nabla_{v_\beta} f
\in C_0(R) \otimes \Omega.
\]
We naturally extend $ \omega $ to $ C_0(R; \mathcal{B}) $
by applying $ \omega $ pointwise.
Since $ E K^{-1} \phi_\alpha \in \ker(\omega) $,
recalling the expression for $ A f $ from \cref{eq:limit-gener-defn}, we get
\[ \begin{split}
\omega(g) &= \omega(L f) - \eta^{-1} \sum_{\alpha=1}^n \omega(E K^{-1} \phi_\alpha) \nabla_{v_\alpha} f
- \eta^{-1} \sum_{\alpha,\beta=1}^n \omega(\phi_\alpha \cdot K^{-1} \phi_\beta) \;\; (\nabla_{v_\alpha} \nabla_{v_\beta} f) \\
&= \omega(L f) + \eta^{-1} \sum_{\alpha,\beta=1}^n \langle \phi_\alpha, \phi_\beta \rangle_K \;\; (\nabla_{v_\alpha} \nabla_{v_\beta} f) = A f \in C_0(R).
\end{split} \]
Since $ 1 \in \Omega $, we have $ \omega(g) \in C_0(R) \subseteq C_0(R) \otimes \Omega $
and so $ \omega(g) - g \in C_0(R) \otimes \Omega $.
Because $ \omega $ is a projection, we have $ \omega(\omega(g) - g) = \omega(g) - \omega(g) = 0 $
so $ \omega(g) - g \in \ker(\omega) \cap (C_0(R) \otimes \Omega) $.
We extend $ K^{-1} $ to $ \ker(\omega) \cap (C_0(R) \otimes \Omega) $ by linearity
and define $ f_2 := \eta^{-1} K^{-1} (\omega(g) - g) \in \ker(\omega) \cap (C_0(R) \otimes \Omega) $.

Now, assume that $ K^{-1} $ is bounded on $ \ker{\omega} \cap \Omega $ and $ L(D) \subseteq C_0(R; \Omega) $.
Then, $ K^{-1} $ can be extended to all of $ \ker(\omega) \cap C_0(R; \Omega) $
and we allow $ g \in C_0(R; \Omega) $.
Then,
\begin{equation} \label{eq:zeroth-ord-val}
g + \eta K f_2 = g + \eta K \eta^{-1} K^{-1} (\omega(g) - g) = g + \omega(g) - g = \omega(g) = A f. 
\end{equation}

Using $ K f = f K(1) = 0 $, \cref{eq:first-ord-vanish}, and \cref{eq:zeroth-ord-val}, it follows that
\[ \begin{split}
\eps{A} \eps{f}
&= \left(L + \frac{1}{\epsilon} (E + J) + \frac{1}{\epsilon^2} \eta K\right)
(f + \epsilon f_1 + \epsilon^2 f_2) \\
&= \frac{1}{\epsilon^2} \eta K f + \frac{1}{\epsilon} \left(\big. \eta K f_1 + (E + J) f \right) \\
&\quad + \left(\big. L f + (E + J) f_1 + \eta K f_2 \right)
+ \epsilon \left(\big. J f_2 + L f_1 + \epsilon L f_2 \right) \\
&= 0 + 0 + (g + \eta K f_2) + O(\epsilon) = A f + O(\epsilon) \to A f
\end{split} \]
as $ \epsilon \downarrow 0 $
proving the condition necessary to apply theorem 6.1 of \cite{EK}.
Thus, for every $ f \in C_0(R) $, $ \eps{T_t} f \to T_t f $ uniformly on bounded intervals of $ t \in [0, \infty) $.
Then, by Corollary 8.7 in chapter 4 of \cite{EK},
since $ C_c^\infty(R) \subseteq C_0(R) $
strongly separates points in $ \mathbb{R}^{d_X} $
(using bump functions),
we conclude that $ \eps{Q} \to Q $ in law.
\end{proof}

Having identified the limiting process in \cref{prop:semigrp-average},
we now present it as the solution to an SDE.
\medskip

\begin{prop} \label{prop:ito-to-strat-DE}
Suppose the process $ Q $ is a Markov process
with infinitesimal generator $ A $
as defined in \cref{eq:limit-gener-defn}, 
$ \eta \in C_b^2(R) $,
and $ \{v_\alpha\}_{\alpha=1}^n \subseteq C_b^3(R; \mathbb{R}^d) $.
Additionally, we assume $ \langle \cdot, \cdot \rangle_K $ is symmetric and positive semi-definite on $ \{\phi_\alpha\} $.
Let $ \omega $ and $ \Omega $ be defined as in \cref{prop:semigrp-average}
and $ L = b \cdot \nabla_R $ for some $ b \in C_b^2(R; \mathbb{R}^d \otimes \Omega) $.
Define $ \wbr{b} := \omega(b) \in C_b(R; \mathbb{R}^d) $.
Then, $ A $ can be written as in \cref{eq:decoupled-vec-fields},
$ D $ is a core for $ A $,
and $ Q $ solves the Stratonovich SDE:
\begin{equation} \label{eq:strat-limit-sde}
\dr Q(t) = \wbr{b}(Q(t)) \dr t
+ \frac{1}{\sqrt{\eta(Q(t))}} \sum_{\alpha=1}^n \sqrt{2} v_\alpha(Q(t)) \circ \dr \tilde{B}_\alpha(t)
\end{equation}
where $ \{\tilde{B}_\alpha\}_{\alpha=1}^n $
is a multivariate Wiener process on $ \mathbb{R}^n $
with covariance $ \expect{\tilde{B}_\alpha(s) \tilde{B}_\beta(t)} = \langle \phi_\alpha, \phi_\beta \rangle_K \cdot \min(s, t) $.
\end{prop}
\begin{proof}
Because $ L = b \cdot \nabla_R $ and $ \omega $ is $ C_b(R) $-linear,
we know $ \omega(L f) = \omega(b \cdot \nabla_R f) = \wbr{b} \cdot \nabla_R f $.
We define the matrix $ C = \{c_{\alpha \beta}\}_{\alpha,\beta=1}^n := \{\langle \phi_\alpha, \phi_\beta \rangle_K\}_{\alpha,\beta=1}^n $.
By assumption, $ C $ is symmetric and positive semi-definite so there exists
a symmetric positive semi-definite matrix $ H = \{h_{\alpha \gamma}\}_{\alpha,\gamma=1}^n := \sqrt{C} $.
Then, since $ C = H^2 $, for all $ \alpha $ and $ \beta $,
we have $ c_{\alpha \beta} = \sum_{\gamma=1}^n h_{\alpha \gamma} h_{\gamma \beta} $.
Using the definition of $ \{\tilde{v}_\gamma\} $, for all $ f \in D $,
\[ \begin{split}
\sum_{\alpha,\beta=1}^n \langle \phi_\alpha, \phi_\beta \rangle_K \nabla_{v_\alpha} \nabla_{v_\beta} f
&= \sum_{\alpha,\beta=1}^n c_{\alpha \beta} \nabla_{v_\alpha} \nabla_{v_\beta} f
= \sum_{\alpha,\beta=1}^n \sum_{\gamma=1}^n h_{\alpha \gamma} h_{\beta \gamma} 
\nabla_{v_\alpha} \nabla_{v_\beta} f \\
&= \sum_{\gamma=1}^n \left(\sum_{\alpha=1}^n \nabla_{h_{\alpha \gamma} v_\alpha} \right)
\left(\sum_{\beta=1}^n \nabla_{h_{\beta \gamma} v_\beta}\right) f
= \sum_{\gamma=1}^n \nabla_{\tilde{v}_\gamma} \nabla_{\tilde{v}_\gamma} f \\
&= \sum_{\gamma=1}^n \left[\Big.
(\nabla_{\tilde{v}_\gamma} \tilde{v}_\gamma) \cdot \nabla_R f
+ (\nabla \nabla_R f)(\tilde{v}_\gamma, \tilde{v}_\gamma)
\right]
\end{split} \]
which, along with the commutation of $ \{\nabla_{v_\alpha}\} $ with $ \eta^{-1/2} $,
demonstrates the equation for $ A f $ in \cref{eq:decoupled-vec-fields}.
Consequently,
\[
A f = \wbr{b} \cdot \nabla_R f + \sum_{\gamma=1}^n \left[
\frac{1}{2} \left(\nabla_{\sqrt{2 \eta^{-1}} \tilde{v}_\gamma} \sqrt{2 \eta^{-1}} \tilde{v}_\gamma\right) \cdot \nabla_R f
+ \frac{1}{2} (\nabla \nabla_R f)\left(\sqrt{2 \eta^{-1}} \tilde{v}_\gamma, \sqrt{2 \eta^{-1}} \tilde{v}_\gamma\right)
\right]
\]
which, by It\^o's formula, is the infinitesimal generator for the It\^o SDE:
\begin{equation} \label{eq:ito-limit-sde}
\dr Q(t) = \wbr{b}(Q(t)) \dr t + \sum_{\gamma=1}^n \left[
\frac{1}{2} \left(\nabla_{\sqrt{2 \eta^{-1}} \tilde{v}_\gamma} \sqrt{2 \eta^{-1}} \tilde{v}_\gamma\right)(Q(t)) \dr t
+ \sqrt{2 \eta^{-1}} \tilde{v}_\gamma(Q(t)) \dr B_\gamma(t)
\right]
\end{equation}
where $ \{B_\gamma\}_{\gamma=1}^n $ is a collection of independent standard Wiener processes
so that $ \expect{B_\gamma(s) B_\gamma(t)} = \min(s, t) $ for all $ \gamma $.
The $ \nabla_{\sqrt{2} \tilde{v}_\gamma} \sqrt{2} \tilde{v}_\gamma $ term
is a Stratonovich correction (see section 4.9 of \cite{KP}).
Additionally, for each $ \alpha $,
we define $ \tilde{B}_\alpha := \sum_{\gamma=1}^n h_{\alpha \gamma} B_\gamma $
so that we can rewrite \cref{eq:ito-limit-sde} as the Stratonovich SDE:
\begin{equation} \label{eq:linear-noise-trf} \begin{split}
\dr Q(t) = \wbr{b}(Q(t)) \dr t
+ \sum_{\gamma=1}^n \sqrt{2 \eta^{-1}} \tilde{v}_\gamma(Q(t)) \circ \dr B_\gamma(t)
&= \wbr{b}(Q(t)) \dr t
+ \frac{1}{\sqrt{\eta(Q(t))}} \sum_{\gamma=1}^n \sum_{\alpha=1}^n \sqrt{2} h_{\alpha \gamma} v_\alpha(Q(t)) \circ \dr B_\gamma(t) \\
&= \wbr{b}(Q(t)) \dr t
+ \frac{1}{\sqrt{\eta(Q(t))}} \sum_{\alpha=1}^n \sqrt{2} v_\alpha(Q(t)) \circ \dr \tilde{B}_\alpha(t)
\end{split} \end{equation}
where $ \{\tilde{B}_\alpha\}_{\alpha=1}^n $
is a multivariate Wiener process
with covariance
\[
\expect{\tilde{B}_\alpha(s) \tilde{B}_\beta(t)}
= \sum_{\gamma=1}^n h_{\alpha \gamma} h_{\gamma \beta} \expect{B_\gamma(s) B_\gamma(t)}
= \min(s, t) \sum_{\gamma=1}^n h_{\alpha \gamma} h_{\gamma \beta}
= c_{\alpha \beta} \cdot \min(s, t)
= \langle \phi_\alpha, \phi_\beta \rangle_K \cdot \min(s, t)
\]
for all $ 1 \leq \alpha, \beta \leq n $ and $ s, t \geq 0 $
proving \cref{eq:strat-limit-sde}.
Lastly, we show that $ D $ is a core for $ A $.
For any $ f \in C_c^\infty(R) $ (i.e. smooth compactly supported functions),
let $ (r, t) \mapsto u(r, t) $ be a solution
to the parabolic PDE $ \partial_t u = T_t u $
with initial condition $ u(r, 0) = f(r) $.
Then, by parabolic regularity,
since $ \wbr{b} $ and $ \eta $ are second-differentiable
and $ \{v_\alpha\} $ are third-differentiable,
for any $ t \geq 0 $, we know $ r \mapsto u(r, t) $ is second-differentiable and vanishes at infinity.
Thus, from proposition 3.3 of \cite{EK}
(with $ D_0 = C_c^\infty(R) $),
we see that $ D $ is a core for the infinitesimal generator $ A $.
\end{proof}

Lastly, in some of the cases, we wish to take $ \mathcal{B} = L^p(S, \mu) $.
Hence, it is necessary to show that the semigroup
for $ \eps{Y} $ is strongly continuous on $ C_0(R, L^p(S, \mu)) $.

\begin{lem} \label{lem:Lp-strong-cont}
Suppose $ Y $ is a Markov process on $ R \times S $ with semigroup $ \{T_t\}_{t \geq 0} $
and $ \mu $ is a $ \{T_t\} $-invariant probability measure on $ S $.
If $ Y $ is a Feller-Dynkin process (i.e. $ t \mapsto T_t $ is strongly continuous as a semigroup on $ C_0(R \times S) $)
then for any $ p \geq 1 $, $ \{T_t\} $ is strongly continuous
as a semigroup on $ C_0(R; L^p(S, \mu)) $.
\end{lem}
\begin{proof}
Note that $ C_c(R) $ (i.e. compactly supported continuous functions)
is dense in $ C_0(R) $
and $ C_c(S) $ is dense in $ L^p(S, \mu) $.
Then, $ C_c(R) \times C_c(S) \subseteq C_c(R \times S) $
is dense in $ C_0(R) \otimes L^p(S, \mu) $
which is, in turn, dense in $ C_0(R; L^p(S, \mu)) $.
Thus, $ C_c(R \times S) $ is dense in $ C_0(R; L^p(S, \mu)) $.
We use $ \|\cdot\|_\infty $ to denote the norm on $ C_0(R \times S) $
and $ \|\cdot\|_p $ for the norm on $ C_0(R; L^p(S, \mu)) $.
\smallskip

Because $ \{T_t\} $ is a semigroup, it suffices to prove strong continuity at $ t = 0 $.
Let $ f \in C_0(R; L^p(S, \mu)) $ and $ \epsilon > 0 $.
By density, there exists $ g \in C_c(R \times S) $
such that $ \|f - g\|_p \leq \epsilon / 3 $.
Then, since $ \{T_t\} $ is strongly continuous on $ C_0(R \times S) $,
there exists $ \delta > 0 $ such that $ \|T_t g - g\|_\infty < \epsilon / 3 $ for all $ 0 \leq t < \delta $.
For any $ f \in L^p(S, \mu) $, since $ \mu $ is an invariant measure, using Jensen's inequality,
\[
\|T_t f\|_p^p = \int_S \left|\int_S f(y) \cdot p_t(x, \!\dr y)\right|^p \mu(\!\dr x)
\leq \int_S |f(y)|^p \int_S p_t(x, \!\dr y) \mu(\!\dr x)
= \int_S |f(y)|^p \cdot \mu(\!\dr y) = \|f\|_p^p
\]
so $ \{T_t\} $ is contractive on $ L^p(S, \mu) $.
Thus, for any $ t \in [0, \delta) $,
\[
\|T_t f - f\|_p
\leq \|T_t f - T_t g\|_p + \|T_t g - g\|_p + \|g - f\|_p
\leq \|T_t(f - g)\|_p + \|T_t g - g\|_\infty + \|g - f\|_p
\leq \frac{\epsilon}{3} + \frac{\epsilon}{3} + \frac{\epsilon}{3} = \epsilon.
\]
Therefore, $ t \mapsto T_t f $ is continuous
and so $ \{T_t\} $ is strongly continuous with respect to $ \|\cdot\|_p $.
\end{proof}

\subsection{Amplitude Scaling}

The following lemma is elementary.  We provide its proof for completeness.

\begin{lem} \label{lem:periodic-diff}
Let $ \phi : \mathbb{R} \to \mathbb{C} $ be continuous and $ 2\pi $-periodic.
If $ \int_0^{2\pi} \phi(\theta) \dr\theta = 0 $
then there exist $ 2\pi $-periodic continuous functions $ \Psi, \Phi : \mathbb{R} \to \mathbb{C} $
such that $ \partial_\theta^2 \Psi = \partial_\theta \Phi = \phi $, $ \|\Phi\| \leq  2\pi \|\phi\| $, $ \|\Psi\| \leq  4\pi^2\|\phi\| $ and
\[ \int_0^{2\pi} \Psi(\theta) \dr\theta = \int_0^{2\pi} \Phi(\theta) \dr\theta = 0. \]
\end{lem}
\begin{proof}
For $\theta \in [0,2\pi]$, we define 
\[
\Phi(\theta) := \int_0^{2\pi} k(\theta, \alpha) \phi(\alpha) \dr\alpha
\quad \text{where} \quad
k(\theta, \alpha) := \frac{1}{2\pi} \begin{cases}
\alpha & \text{for } \alpha \leq \theta \\
\alpha - 2\pi & \text{for } \alpha > \theta
\end{cases}.
\]
Since $ \Phi(0) = \Phi(2\pi) $, $ \Phi $ can be extended to a continuous $ 2\pi $-periodic function on $ \mathbb{R} $.
Morally, $ (\partial_\theta k)(\theta, \alpha) = \delta(\theta - \alpha) $ so $ \Phi' = \phi $.
Since $ |k(\theta, \alpha)| \leq 1 $, we have $ \|\Phi\| \leq 2\pi \|\phi\| $.
Because $ \int_0^{2\pi} k(\theta, \alpha) \dr\theta = 0 $, we see that $ \Phi $ is mean zero.
Applying the same construction to $ \Phi $ gives $ \Psi $.
\end{proof}

\begin{proof}[Proof of \cref{thm:ampl-scaling-limit}]
We prove the convergence in law of the processes using \cref{prop:semigrp-average}.
First, we introduce the random process $ \eps{\Theta}(t) = \frac{1}{\epsilon} \vartheta(M(t)) $
so that \cref{eq:ampl-scaling} becomes
\begin{equation} \label{eq:ampl-scaling-extended} \begin{split}
\dr \eps{X}(t) &= b\left( \eps{X}(t), \eps{\Theta}(t) \right) \dr t
+ \frac{1}{\epsilon} \sum_{\alpha=1}^{d_B} v_\alpha\left( \eps{X}(t) \right) \phi_\alpha\left( \eps{\Theta}(t) \right) \dr t \\
\dr M(t) &= \nu(M(t)) \dr t + \sum_{k=1}^{d_W} \sigma_k(M(t)) \dr W_k(t) \\
\dr \eps{\Theta}(t) &= \frac{1}{\epsilon} \rho(M(t)) \dr t + \frac{1}{\epsilon} \sum_{k=1}^{d_W} \varsigma_k(M(t)) \dr W_k(t).
\end{split} \end{equation}
where the functions $ \rho : \mathbb{R}^{d_M} \to \mathbb{R} $
and $ \{\varsigma_k : \mathbb{R}^{d_M} \to \mathbb{R}\}_{k=1}^{d_W} $
are defined as
\[ \begin{split}
\rho(m) &:= \nu(m) \cdot \nabla \vartheta(m)
	+ \frac{1}{2} \sum_{k=1}^{d_W} \sum_{i,j=1}^{d_M} \sigma_{k i}(m) \sigma_{k j}(m) \;\;
    (\partial_{m_i} \partial_{m_j} \vartheta)(m) \\
\varsigma_k(m) &:= \sigma_k(m) \cdot \nabla \vartheta(m)
\quad \text{ for all } k \in \{1,\dots, d_W\}.
\end{split} \]

Note that $ M $ is a Markov process.
While neither $ \eps{X}(t) $ nor $ \eps{\Theta}(t) $
are Markovian on their own, the whole system $ \eps{Y}(t) := \left(\eps{X}(t), M(t), \eps{\Theta}(t)\right) $ is Markovian.
Since only $ 2\pi $-periodic functions of $ \eps{\Theta}(t) $ are considered,
we can take the state space of $ \eps{Y}(t) $ to be $ \mathbb{R}^{d_X} \times \mathbb{R}^{d_M} \times \mathbf{S}^1 $
where we identify $ \mathbf{S}^1 $ with $ \mathbb{R} / 2\pi \mathbb{Z} $.
We represent points of this space using the variables $ (x, m, \theta) $.
$ \eps{Y}(t) $ is a Feller-Dynkin process in the sense of \cite{RW}.
Let $ T^{(\epsilon)} $ be the associated semigroup on the Banach space $ C_0(\mathbb{R}^{d_X} \times \mathbb{R}^{d_M}; C_b(\mathbf{S}^1)) = C_0(\mathbb{R}^{d_X} \times \mathbb{R}^{d_M} \times \mathbf{S}^1) $.
We denote its infinitesimal generator by $ \eps{A} $.
For differentiable functions $ f $, we write $ \nabla_x f = (\partial_{x_1}, \ldots, \partial_{x_{d_X}}) $
and similarly for $ \nabla_m f $.
It follows from the It\^o formula that for sufficiently differentiable $ f $,
\[ \begin{split}
\left(\eps{A} f \right)(x, m, \theta) &=
b(x, \theta) \cdot \nabla_x f + \frac{1}{\epsilon} \sum_{\alpha=1}^{d_B} \phi_\alpha(\theta) v_\alpha(x) \cdot \nabla_x f
+ \frac{1}{\epsilon} \rho(m) \partial_\theta f
+ \nu(m) \cdot \nabla_m f \\
&\qquad + \frac{1}{2} \sum_{k=1}^{d_W} \left[
	\frac{1}{\epsilon^2} \varsigma_k(m)^2 \partial_\theta^2 f
	+ \frac{2}{\epsilon} \varsigma_k(m) (\sigma_k(m) \cdot \nabla_m (\partial_\theta f))
	+ \sum_{i,j=1}^{d_M} \sigma_{k i}(m) \sigma_{k j}(m) \;\; \partial_{m_i} \partial_{m_j} f
\right].
\end{split} \]
The process $ M $ is a Feller-Dynkin process in its own right
and so has an infinitesimal generator $ A_M $.
For sufficiently differentiable $ f $, we have
\[
(A_M f)(x, m, \theta) = \nu(m) \cdot \partial_m f + \frac{1}{2} \sum_{k=1}^{d_W}
\sum_{i,j=1}^{d_M} \sigma_{k i}(m) \sigma_{k j}(m) \;\;\; \partial_{m_i} \partial_{m_j} f.
\]
We wish to apply \cref{prop:semigrp-average} so we decompose $ \eps{A} $ as
\[ \begin{split}
\eps{A} f = L f + \frac{1}{\epsilon} (J &+ E) f + \frac{1}{\epsilon^2} \eta K f
\quad \text{for sufficiently differentiable $ f $ where} \\
L := b(x,\theta) \cdot \nabla_x + A_M
&\qquad
J := \sum_{\alpha=1}^{d_B} \phi_\alpha(\theta) v_\alpha(x) \cdot \nabla_x \\
E := \rho(m) \partial_\theta
+ \sum_{k=1}^{d_W} \varsigma_k(m) (\sigma_k(m) \cdot \nabla_m) \partial_\theta
&\qquad
\eta(m) := \frac{1}{2} \kappa(m)^2 = \frac{1}{2} \sum_{k=1}^{d_W} \varsigma_k(m)^2
\qquad
K := \partial_\theta^2
\end{split} \]
Now, we take $ R = \mathbb{R}^{d_X} \times \mathbb{R}^{d_M} $ and $ S = \mathbf{S}^1 $.
Also, we take $ \mu $ to be the uniform probability measure on $ \mathbf{S}^1 $
and $ \mathcal{B} := \Omega := C_b(S) $.

Next, we verify that the operators fulfill the necessary conditions.
\begin{enumerate}
\item $ \partial_\theta^2 (1) = 0 $
and, by \cref{lem:periodic-diff}, $ \partial_\theta^2 $ is invertible on $ \ker(\omega) \cap \Omega = \ker(\omega) \cap C_b(S) $.
\item Since $ \eta^{-1/2} = \sqrt{2} / \kappa $ is a function of only $ m $, $ \nabla_x \eta^{-1/2} = 0 $.
Since $ \phi_\alpha, K^{-1} \phi_\beta \in C_b(S) $,
we have $ \phi_\alpha \cdot K^{-1} \phi_\beta \in C_b(S) = \Omega $.
\item $ \partial_\theta $ is $ C_0(R) $-linear so $ E $ is as well and, further, $ C_0(R) \subseteq \ker(\partial_\theta) \subseteq \ker(E) $.
If $ \int \phi \dr\theta = 0 $ and $ \phi $ is differentiable
then $ \int \partial_\theta \phi \dr\theta = 0 $ too
so $ E $ preserves $ \ker(\omega) \cap \Omega $.
\item The elements of $ D $ have first derivatives in $ C_0(R) $ so $ D \subseteq \domain{L} $.
Again, by \cref{lem:periodic-diff}, $ K^{-1} $ is bounded on $ \ker(\omega) \cap \Omega $
and, by the assumptions on $ b $, $ L f = b \cdot \nabla_x f \in C_0(R \times S) = C_0(R;\Omega) $.
Since $ \nabla_x $ is $ \mathcal{H} $-linear, so is $ L $.
\end{enumerate}

For each $ \alpha $, let $ \Psi_\alpha, \Phi_\alpha \in C_b(S) $
such that $ \partial_\theta^2 \Psi_\alpha = \partial_\theta \Phi_\alpha = \phi_\alpha $.
Then, the bilinear form $ \langle \cdot, \cdot \rangle_K $ is
\[ \begin{split}
c_{\alpha \beta} = \langle \phi_\alpha, \phi_\beta \rangle_K
&= -\mu(\phi_\alpha \cdot K^{-1} \phi_\beta)
= - \frac{1}{2\pi} \int_0^{2\pi} \phi_\alpha(\theta) \cdot \Psi_\beta(\theta) \dr\theta \\
&= - \frac{1}{2\pi} \phi_\alpha(\theta) \Phi_\beta(\theta) \Big|_{\theta=0}^{2\pi}
+ \frac{1}{2\pi} \int_0^{2\pi} \Phi_\alpha(\theta) \cdot \Phi_\beta(\theta) \dr\theta
= \frac{1}{2\pi} \langle \Phi_\alpha, \Phi_\beta \rangle_{L^2([0, 2\pi))}
\end{split} \]
By \cref{prop:ito-to-strat-DE},
we know the solution $ X $ to the SDE (\ref{eq:ampl-limit-sde})
has infinitesimal generator $ A $ fulfilling \cref{eq:limit-gener-defn}.
Then, by \cref{prop:semigrp-average}, the solutions $ \eps{X} $ to the SDE (\ref{eq:ampl-scaling})
converge in law to $ X $.
\end{proof}

\subsection{Non-Compact Time Scaling by Ergodicity}

\Cref{lem:ergodic-inv} serves a similar purpose as \cref{lem:periodic-diff}
did for the amplitude scaling argument---it ensures the existence of pre-images.

\begin{proof}[Proof of \cref{lem:ergodic-inv}]
Let $ T_M $ be the semigroup for $ M $ on $ \mathcal{B}_M $.
We define the operator $ A_M^{-1} : D_0 \to D_0 $ as
\[ A_M^{-1} \phi := -\int_0^\infty T_M(t) \phi \dr t \in D_0. \]
for all $ \phi \in D_0 \subseteq \mathcal{B}_M $.
It remains to show that $ A_M^{-1} $ is well-defined
while its linearity is immediate.
For any $ t \geq 0 $ and $ m \in S $, let $ \mu_{t,m} $
be the distribution of $ M(t) $ given that $ M(0) = m $.
Thus, $ \mu_{t,m}(E) = \probWhen{M(t) \in E}{M(0) = m} $ for every $ E \in \mathcal{M} $.
From the strong ergodicity of $ M $,
we get $ \|\mu_{t,m} - \mu\| \leq C e^{-\lambda t} $ for any $ m \in S $ and $ t \geq t_0 $.
For any $ t \geq t_0 $, since $ \int \phi \dr\mu = 0 $,
\begin{equation} \label{eq:semigrp-bnd} \begin{split}
\left|\Big. (T_M(t) \phi)(m) \right|
&= \left|\Big.\expectWhen{\Big. \phi(M(t))}{M(0) = m} \right|
= \left| \int \phi \dr\mu_{t,m} \right| \\
&= \left| \int \phi \dr\mu_{t,m} - \int \phi \dr\mu \right| \\
&\leq \int |\phi| \cdot \dr|\mu_{t,m} - \mu|
\leq \|\phi\|_\infty \cdot C e^{-\lambda t}.
\end{split} \end{equation}
Thus, we can bound the integral with respect to $ t $ by
\[ \begin{split}
\int_0^\infty \left\|\Big. T_M(t) \phi \right\|_\infty \dr t
&= \int_0^{t_0} \left\|\Big. T_M(t) \phi \right\|_\infty \dr t
+ \int_{t_0}^\infty \left\|\Big. T_M(t) \phi \right\|_\infty \dr t \\
&\leq t_0 \|\phi\|_\infty
+ \|\phi\|_\infty \int_0^\infty C e^{-\lambda t} \dr t < \infty
\end{split} \]
implying that $ A_M^{-1} \phi $ is well-defined.
From the above inequality, for any $ \phi \in D_0 $, we also have
\[
\|A_M^{-1} \phi\|_\infty \leq t_0 \|\phi\|_\infty
+ \|\phi\|_\infty \int_0^\infty C e^{-\lambda t} \dr t
= \left(t_0 + \frac{C}{\lambda}\right) \|\phi\|_\infty
\]
so $ \|A_M^{-1}\| \leq t_0 + \frac{C}{\lambda} < \infty $
and $ A_M^{-1} $ is bounded.
Also, using the invariance of $ \mu $, we get
\[
\int_S (A_M^{-1} \phi)(m) \dr\mu(m)
= -\int_0^\infty \int_S T_M(t) \phi \dr\mu \dr t
= \int_0^\infty \left(\int_S \phi \dr\mu\right) \dr t = 0
\]
so the range of $ A_M^{-1} $ lies in $ D_0 $.
Lastly, we know $ t \mapsto T_M(t) \phi $ is continuous so
\[ \begin{split}
A_M A_M^{-1} \phi
&= -\lim_{h \downarrow 0} \frac{T_M(h) - I}{h}
\int_0^\infty T_M(t) \phi \dr t
= \lim_{h \downarrow 0} \frac{1}{h} \left[
\int_0^\infty T_M(t) \phi \dr t
- \int_0^\infty T_M(t+h) \phi \dr t
\right] \\
&= \lim_{h \downarrow 0} \frac{1}{h} \left[
\int_0^\infty T_M(t) \phi \dr t
- \int_h^\infty T_M(t) \phi \dr t
\right]
= \lim_{h \downarrow 0} \frac{1}{h} \int_0^h T_M(t) \phi \dr t
= T_M(0) \phi = \phi
\end{split} \]
implying that $ A_M A_M^{-1} = I_{D_0} $.
\end{proof}
\bigskip

\begin{proof}[Proof of \cref{thm:ergodic-time-scaling}]
We can rewrite \cref{eq:time-scaling} by defining $ \eps{M}(t) = M(t / \epsilon^2) $
to obtain
\begin{equation} \label{eq:time-scaling-mod}
\dr \eps{X}(t) = b\left(\eps{X}(t), \eps{M}(t)\right) \dr t
+ \frac{1}{\epsilon} \sum_{\alpha=1}^{d_B} \phi_\alpha\left(\eps{M}(t)\right)
\, v_\alpha\left( \eps{X}(t) \right) \dr t.
\end{equation}
By assumption, $ M $ is a Feller-Dynkin process
so let $ \{T_M(t)\}_{t \geq 0} $
be its Feller semigroup on $ \mathcal{B}_M \subseteq C_b(S) $
and $ A_M $ be the associated infinitesimal generator.
If we take $ T_{\eps{M}} $ to be the semigroup of $ \eps{M} $ then for any $ m \in S $,
\[
(T_{\eps{M}}(t) f)(m) = \expectWhen{f\left( \eps{M}(t) \right)}{\eps{M}(0) = m}
= \expectWhen{f(M(t / \epsilon^2))}{M(0) = m}
= (T_M(t / \epsilon^2) f)(m).
\]
Therefore, the infinitesimal generator $ A_{\eps{M}} $ of $ \eps{M} $ is
\[
A_{\eps{M}} f = \lim_{t \downarrow 0} \frac{T_{\eps{M}}(t) f - f}{t}
= \lim_{t \downarrow 0} \frac{T_M(t / \epsilon^2) f - f}{t}
= \lim_{t \downarrow 0} \frac{T_M(t) f - f}{\epsilon^2 t}
= \frac{1}{\epsilon^2} A_M f
\]
for all $ f \in \domain{A_M} $.
Now, $ \left( \eps{X}, \eps{M} \right) $ is a Feller-Dynkin process with infinitesimal generator $ \eps{A} $
which satisfies
\[
\eps{A} f = b(x, m) \cdot \nabla_x f
+ \frac{1}{\epsilon} \sum_{\alpha=1}^{d_B} \phi_\alpha(m) v_\alpha(x) \cdot \nabla_x f
+ \frac{1}{\epsilon^2} A_M f
\]
for all sufficiently differentiable $ f $.
We wish to apply \cref{prop:semigrp-average}.
We take $ R = \mathbb{R}^{d_X} $
and $ S $ is given.
We take $ \mu $ to be the invariant measure for $ M $
and $ \mathcal{B} := \Omega := \mathcal{B}_M $.
Next, we decompose $ \eps{A} $ for sufficiently differentiable $ f $ as
\[ \begin{split}
\eps{A} f = L f + \frac{1}{\epsilon} &J + \frac{1}{\epsilon^2} K f
\quad \text{where} \\
L := b(x, m) \cdot \nabla_x
\qquad
J := \sum_{\alpha=1}^{d_B} &\phi_\alpha(m) v_\alpha(x) \cdot \nabla_x
\qquad
K := A_M
\end{split} \]
while $ \eta = 1 $ and $ E = 0 $.
We now verify that these operators fulfill the conditions of \cref{prop:semigrp-average}.
\begin{enumerate}
\item $ K(1) = A_M(1) = 0 $ and, by \cref{lem:ergodic-inv},
$ A_M $ is right-invertible on $ \ker(\omega) \cap \Omega $.
\item For all $ \alpha $ and $ \beta $, $ \nabla_{v_\alpha} (\eta^{-1/2}) = \nabla_{v_\alpha} (1) = 0 $
while $ \phi_\alpha, K^{-1} \phi_\beta \in \mathcal{B}_M $
which, since $ \mathcal{B}_M $ is an algebra, means $ \phi_\alpha \cdot K^{-1} \phi_\beta \in \mathcal{B}_M = \Omega $.
\item Since $ E = 0 $, it fulfills all of the necessary properties.
\item The elements of $ D $ have first derivatives in $ C_0(R) $ so $ D \subseteq \domain{L} $.
From \cref{lem:ergodic-inv}, we know $ K^{-1} = A_M^{-1} $ is bounded on $ \ker(\omega) \cap \Omega $
and, by the assumptions on $ b $, $ L f = b \cdot \nabla_x f \in C_0(R; \mathcal{B}_M) = C_0(R; \Omega) $.
Since $ \nabla_x $ is $ \mathcal{H} $-linear so is $ L $.
\end{enumerate}
By \cref{prop:ito-to-strat-DE}, we see that the diffusion process solving $ X $ of the SDE (\ref{eq:time-limit-sde})
has an infinitesimal generator $ A $ fulfilling \cref{eq:limit-gener-defn}.
Thus, by \cref{prop:semigrp-average}, the solutions $ \eps{X} $ to the SDE (\ref{eq:time-scaling-mod})
converge in law to $ X $.
\end{proof}
\bigskip

\begin{proof}[Proof of \cref{coro:ergodic-wiener-well}]
Let $ A_M $ be the infinitesimal generator for $ M $ so that
\[
A_M f = \frac{1}{2} \nabla_m^2 f - (\nabla_m U) \cdot (\nabla_m f)
\]
for all sufficiently differentiable $ f \in C_b(\mathbb{R}^{d_M}) $.
By 1.3B of \cite{Dynk}, $ \mathcal{B}_M = \overline{\domain{A_M}} $
so $ \{\phi_\alpha\} \subseteq \domain{A_M} \subseteq \mathcal{B}_M $.
Further, second differentiable functions on $ \mathbb{R}^{d_M} $ 
with appropriately bounded derivatives form a core for $ A_M $
and, since these functions form an algebra, $ \mathcal{B}_M $ is an algebra as well.
Then, we can use $ C_0^2(\mathbb{R}^{d_M}) $ as a core
and for any $ f \in C_0^2(\mathbb{R}^{d_M}) $,
\[
\int (A_M f)(m) \frac{e^{-2 U(m)}}{Z} \dr m
= \int \left(\frac{1}{2} \nabla^2 f - (\nabla U) \cdot (\nabla f) \right) \frac{e^{-2 U}}{Z} \dr m
= \int f \left(\frac{1}{2} \nabla^2 e^{-2 U} + \nabla \cdot ((\nabla U) e^{-2 U}) \right) \frac{1}{Z} \dr m
\]
where $ Z > 0 $ is an appropriate normalization
to ensure $ \frac{1}{Z} e^{-2U} \dr m $ is a probability measure.
Now, we observe that $ \frac{1}{2}\nabla e^{-2U} + \left(\nabla U\right) e^{-2U} = 0 $
which implies that $ \frac{1}{Z} e^{-2U} \dr m $ is an invariant measure for $ M $.
Conjugating $ A_M $
by the multiplication operator $ e^U $, we obtain
\[ \begin{split}
e^{-U} A_M (e^U f) &= \frac{1}{2} e^{-U} \nabla^2 (e^U f) - e^{-U} (\nabla U) \cdot \nabla (e^U f) \\
&= \frac{1}{2} e^{-U} \left[e^U \nabla^2 f + 2 e^U (\nabla U) \cdot (\nabla f) + e^U (|\nabla U|^2 + \nabla^2 U) f \right]
- e^{-U} (\nabla U) \cdot \left[ e^U \nabla f + e^U (\nabla U) f \right] \\
&= \frac{1}{2} \nabla^2 f + (\nabla U) \cdot (\nabla f)
+ \frac{1}{2} (|\nabla U|^2 + \nabla^2 U) f
- (\nabla U) \cdot (\nabla f) - |\nabla U|^2 f \\
&= \frac{1}{2} \nabla^2 f + \frac{1}{2} (\nabla^2 U - |\nabla U|^2) f
= \frac{1}{2} \nabla^2 f - V f
= - H f
\end{split} \]
implying that $ A_M = - e^U H e^{-U} $ and $ A_M^{-1} = - e^U H^{-1} e^{-U} $
where $ H^{-1} $ is the inverse of $ H $ for $ \phi \in D_0 := \{\psi \in \mathcal{B}_M \,:\, 0 = \int_S \psi \dr \mu = \int Z^{-1} e^{-2U} \psi \dr m \}$.
Then, we know that the bilinear form associated with $ M $ is
\[ \begin{split}
\langle \phi, \psi \rangle_M
&= - \int_{\mathbb{R}^{d_M}} \phi(m) (A_M^{-1} \psi)(m) \dr \mu(m)
= \int \phi(m) (e^U H^{-1} e^{-U} \psi)(m) \frac{1}{Z} e^{-2U} \dr m \\
&= \frac{1}{Z} \int (e^{-U} \phi)(m) \cdot (H^{-1} e^{-U} \psi)(m) \dr m.
\end{split} \]
Since $ H $ is self-adjoint with respect to$ \dr m $, the form
$ \langle \cdot, \cdot \rangle_M $ symmetric.
It is known that $ \phi = e^{-U} $ is a non-degenerate eigenvector of $ H $ when the potential $ V = \frac{1}{2} (|\nabla U|^2 - \nabla^2 U) \to \infty $ as $ z \to \infty $ (theorem XIII.47 of \cite{RS4}).
Thus, $ H^{-1} $ is positive on  $ D_0 $ 
and so $ \langle \cdot, \cdot \rangle_M $ is positive semi-definite on $ \{\phi_\alpha\} $.

From the first part of \cref{eq:wang-conds}, for all $ x \in \mathbb{R}^{d_M} $, we get the bound
\begin{equation} \label{eq:dot-prod-bdd}
x \cdot \nabla U(x)
= x \cdot \int_0^1 |x| (\nabla \nabla U(tx)) x \dr t
\geq \int_0^1 |x|^3 h_1'(t |x|) \dr t
= |x|^2 (h_1(|x|) - h_1(0)) = |x|^2 h_1(|x|)
\end{equation}
and similarly $ x \cdot \nabla U(x) \leq |x|^2 h_2(|x|) $.
By theorem 1.1(b) in Wang \cite{JW},
from \cref{eq:dot-prod-bdd}, we know $ h_1 \leq g \leq h_2 $
and so the conditions in \cref{eq:wang-conds}
imply that $ M $ is strongly ergodic.
Thus, from \cref{thm:ergodic-time-scaling},
we obtain the weak convergence.
\end{proof}

\subsection{Time Scaling driven by Integrated Noise}

\begin{lem} \label{lem:schrod-preimg}
Suppose $ \rho : \mathbb{R}^{d_Z} \to \mathbb{R} $ and $ U : \mathbb{R}^{d_Z} \to \mathbb{R} $ are continuous functions
and $ U $ fulfills \cref{eq:poten-conds}.
We assume that there exist $ a, b > 0 $
such that $ V > -b $ and $ |\rho|^2 \leq a (V + b) $.
Let $ f \in \mathcal{H} := L^2(\mathbb{R}^{d_Z}, e^{-2U(z)} \dr z) $.
Then, for any integer $ n \neq 0 $, 
there exists $ g \in \mathcal{H} $ such that
\[
\left[ \big. \frac{1}{2} \nabla_z^2 - (\nabla U) \cdot \nabla_z + \rho \partial_\theta \right] (g e^{n\theta i}) = f e^{n\theta i}
\quad \text{and} \quad \mathrm{Re} \langle f, g \rangle_\mathcal{H} \leq 0.
\]
If in addition $ \langle f, 1 \rangle_\mathcal{H} = 0 $, then there exists a $ g \in L^2(\mathbb{R}^{d_Z}, e^{-2U(z)} \dr z) $ such that 
\[
\left[ \big. \frac{1}{2} \nabla_z^2 - (\nabla U) \cdot \nabla_z + \rho \partial_\theta \right] g = f
\quad \text{and} \quad \mathrm{Re} \langle f, g \rangle_\mathcal{H} \leq 0.
\]
\end{lem}
\begin{proof}
For any $n$, we want to solve the equation 
\[ 
\left[ \big. \frac{1}{2} \nabla_z^2 - (\nabla U) \cdot \nabla_z\right] g + in\rho g = f.
\]
The operator 
\[ g \mapsto Lg = \left[ \big. \frac{1}{2} \nabla_z^2 - (\nabla U) \cdot \nabla_z + \rho \partial_\theta \right] g \]
is a generator of a diffusion in the space $ \mathcal{H} $. In particular, it has no positive spectrum, and thus the same is true for its image under the conjugation by $e^{-U}$, which is a unitary operator between the two Hilbert spaces.  Explicitly, substituting $g = e^U \phi$ and multiplying the resulting operator by $-1$ for greater clarity, we obtain the Schr\"odinger operator $ H $ on the space $ L^2(\mathbb{R}^{d_Z}) $:
\[
\phi \mapsto H\phi = -e^{-U}L\left(e^U \phi\right) = - \frac{1}{2} \nabla_z^2\phi + V\phi
\quad \text{with} \quad
V = \frac{1}{2}\left(\left|\nabla_z U\right|^2 - \nabla_z^2U\right).
\]
Under our assumptions on $U$, the potential $V$ is bounded below and goes to $+\infty$ for $|z| \to \infty$.
It follows that $H$ has only discrete spectrum (theorem 10.4.3 in \cite{Stern}).
An explicit calculation shows that $ H $ applied to the function $ e^{-U} $ is zero, so zero is its eigenvalue.
Since $ H $ has the same spectrum as $ -L $,
it follows that zero is its lowest eigenvalue which by theorem XIII.47 of \cite{RS4} is non-degenerate.
In terms of $ \phi $, the equation for $ g $ becomes
\[H\phi +in\rho\phi = -fe^{-U}\]
If $n = 0$  this equation has a solution if (and only if) its right-hand side is orthogonal to the ground state, i.e. if $ \langle f, 1 \rangle_\mathcal{H} = 0 $.  To treat the case when $n \not = 0$, we will first show that the spectrum of the Schr\"odinger operator with the complex potential $V + in\rho$ is also discrete.  Indeed, denoting the operators of multiplication by $\rho$ and $V$ also by $\rho$ and $V$, we have
\begin{equation} \label{eq:resolv-factors}
\rho\left(bI + H\right)^{-1} = \rho (bI + V)^{-\frac{1}{2}} (bI + V)^{\frac{1}{2}} \left(bI + H\right)^{-\frac{1}{2}} \left(bI + H\right)^{-\frac{1}{2}}
\end{equation}
The product $\rho (bI + V)^{-\frac{1}{2}} $ is bounded by assumption.
The product $ (bI + V)^{\frac{1}{2}} (bI + H)^{-\frac{1}{2}} $ is bounded because $ bI + V \leq bI + H $,
so for any $ \phi $
\[ \begin{split}
\left((bI + V)^{\frac{1}{2}} (bI + H)^{-\frac{1}{2}} \phi, (bI + V)^{\frac{1}{2}} (bI + H)^{-\frac{1}{2}} \phi\right)
&= \left((bI + H)^{-\frac{1}{2}}\phi, (bI + V) (bI + H)^{-\frac{1}{2}} \phi\right) \\
&\leq \left((bI + H)^{-\frac{1}{2}} \phi, (bI + H) (bI + H)^{-\frac{1}{2}} \phi\right) = \left(\phi, \phi\right)
\end{split} \]
Finally, since $ b > 0 $ and the spectrum of $ H $ is non-negative and discrete,
the factor $ (bI + H)^{-\frac{1}{2}} $ in \cref{eq:resolv-factors} is a compact operator.
So the operator $\rho (bI + H)^{-1} $ is a product of a bounded operator and a compact operator,
hence is also compact.
This shows that the perturbation $ in\rho $ is compact relative to $ H $
and thus, by the Weyl theorem (theorem 9.5.1 in \cite{Stern}),
adding it does not change the essential spectrum of $H$.
We claim that $0$ is {\it not} an eigenvalue of the perturbed operator.
For, suppose that $ \psi = \alpha + i\beta $
and $ H\psi +\rho \psi= 0 $.
Taking the real and imaginary parts, we obtain:
\[
H\alpha -\rho\beta = 0
\quad \text{ and } \quad
H\beta + \rho\alpha = 0.
\]
Multiplying the first equation by $ \alpha $ and the second equation by $ \beta $, adding up and integrating, we get
\[
(\alpha, H\alpha) + (\beta, H\beta) = 0
\]
Since $H$ is nonnegative-definite, this implies that $ (\alpha, H\alpha) = (\beta, H\beta) = 0$.  Since the zero eigenvalue of $H$ is simple, this is only possible if  $\alpha$ and $\beta$, and hence also $\psi$, are multiples of the ground state of $H$---a contradiction.  It now follows from the Fredholm alternative that the equation
\[H\phi + in\rho\phi = -fe^{-U} \] has a (unique) solution.
Then, since $ \rho $ is real and $ H $ is positive semi-definite,
\[
\mathrm{Re} \langle f, g \rangle_\mathcal{H}
= \mathrm{Re} \langle -(H + in\rho) \phi, \phi \rangle
= \mathrm{Re} \left(\langle -H \phi, \phi \rangle + in \langle \rho \phi, \phi \rangle\right)
= \langle -H \phi, \phi \rangle \leq 0
\]
where $ \langle \cdot, \cdot \rangle $ is the inner product in $ L^2(\mathbb{R}^{d_Z}) $.
\end{proof}

\begin{proof}[Proof of \cref{thm:integ-noise-time-scaling}]
Let $ \eps{Y} := (\eps{X}, \eps{\Theta}, \eps{Z}) $ be a solution of \cref{eq:integ-noise-time-scaling}.
We regard $ \eps{Y} $ as a process in $ \mathbb{R}^{d_X} \times \mathbf{S}^1 \times \mathbb{R}^{d_Z} $.
Then, $ \eps{Y} $ is a Feller-Dynkin process with generator $ \eps{A} $ given by
\[
\left(\eps{A} f\right)(x, \theta, z) = b(x, \theta) \cdot \nabla_x f
+ \frac{1}{\epsilon} \sum_{\alpha=1}^{d_B} \phi_\alpha(\theta) v_\alpha(x) \cdot \nabla_x f
+ \frac{1}{\epsilon^2} \rho(z) \partial_\theta f
- \frac{1}{\epsilon^2} \nabla U(z) \cdot \nabla_z f + \frac{1}{2 \epsilon^2} \nabla_z^2 f
\]
for sufficiently differentiable $ f \in C_0(\mathbb{R}^{d_X} \times \mathbf{S}^1 \times \mathbb{R}^{d_Z}) $.
We wish to apply \cref{prop:semigrp-average}.
We take $ R := \mathbb{R}^{d_X} $ and $ S := \mathbf{S}^1 \times \mathbb{R}^{d_Z} $.
We define $ \mu_Z $ to be the probability measure on $ \mathbb{R}^{d_Z} $
obtained by normalizing $ e^{-2U(z)} \dr z $.
Then, we define the probability measure $ \mu := \frac{1}{2\pi} \dr \theta \times \mu_Z $
on $ S = \mathbf{S}^1 \times \mathbb{R}^{d_Z} $
and take $ \mathcal{B} := L^2(S, \mu) $ as well as
\[ \Omega := \mathrm{Trig}(\mathbf{S}^1) \otimes L^2(\mathbb{R}^{d_Z}, \mu_Z)
\quad \text{where} \quad
\mathrm{Trig}(\mathbf{S}^1) := \left\{\phi \in C_b(\mathbf{S}^1) \,:\,
\phi \text{ is a trigonometric polynomial}
\right\}.
\]
By \cref{lem:Lp-strong-cont} for $ p = 2 $, the semigroup for $ \eps{Y} $
is strongly continuous on $ C_0(R; \mathcal{B}) $.
Next, we decompose $ \eps{A} $ as
\[ \begin{split}
\eps{A} f = L f + \frac{1}{\epsilon} &J f + \frac{1}{\epsilon^2} K f
\quad \text{for sufficiently differentiable $ f $ where} \\
L := b(x, \theta) \cdot \nabla_x
\qquad
J := \sum_{\alpha=1}^{d_B} &\phi_\alpha(\theta) v_\alpha(x) \cdot \nabla_x
\qquad
K := \frac{1}{2} \nabla_z^2 - \nabla U(z) \cdot \nabla_z + \rho(z) \partial_\theta
\end{split} \]
while $ \eta = 1 $ and $ E = 0 $.
We also define the bilinear form $ \langle \cdot, \cdot \rangle_{\Theta,Z} $
for any $ \phi, \psi \in \ker(\omega) \cap \Omega $ as
\begin{equation} \label{eq:integ-noise-form}
\langle \phi, \psi \rangle_{\Theta,Z} := -\omega(\phi^* \cdot K^{-1} \psi) = \langle \phi, \psi \rangle_K.
\end{equation}
Note that, since $ H $ in the proof of \cref{lem:schrod-preimg} is self-adjoint,
the form $ \langle \cdot, \cdot \rangle_{\Theta,Z} $ is symmetric.
Moreover, because of the negativity conditions in \cref{lem:schrod-preimg},
the form $ \langle \cdot, \cdot \rangle_{\Theta,Z} $ is positive semi-definite.
Now, we verify the conditions for \cref{prop:semigrp-average}.
\begin{enumerate}
\item Observe that $ K(1) = 0 $
and, by \cref{lem:schrod-preimg}, since the elements $ f e^{n\theta i} $
for $ f \in L^2(\mathbb{R}^{d_Z}, \mu_Z) $ span $ \Omega $,
$ K $ is right-invertible on $ \Omega $.
\item For all $ \alpha $ and $ \beta $, $ \nabla_{v_\alpha} (\eta^{-1/2}) = \nabla_{v_\alpha} (1) = 0 $.
Since $ \mathrm{Trig}(\mathbf{S}^1) $ is an algebra
with $ \phi_\alpha \in \mathrm{Trig}(\mathbf{S}^1) \otimes 1 $
and $ K^{-1} \phi_\beta \in \mathrm{Trig}(\mathbf{S}^1) \otimes L^2(\mu_Z) $,
it follows that $ \phi_\alpha \cdot K^{-1} \phi_\beta \in \mathrm{Trig}(\mathbf{S}^1) \otimes L^2(\mu_Z) = \Omega $.
\item Since $ E = 0 $, it fulfills all of the necessary conditions.
\item The elements of $ D $ have first derivatives in $ C_0(R) $ so $ D \subseteq \domain{L} $.
By the assumptions on $ b $, $ L f = b \cdot \nabla_x f \in C_0(R) \otimes \Omega $.
Since $ \nabla_x $ is $ L^2(S, \mu) $-linear so is $ L $.
\end{enumerate}
Then, by \cref{prop:ito-to-strat-DE}, we see that the solution $ X $ of the SDE (\ref{eq:integ-noise-limit-sde})
has an infinitesimal generator $ A $ fulfilling \cref{eq:limit-gener-defn}.
Thus, by \cref{prop:semigrp-average}, the solutions $ \eps{X} $ to the SDE (\ref{eq:integ-noise-time-scaling})
converge in law to $ X $.
\end{proof}

\begin{proof}[Proof of \cref{prop:ou-proc-form}]
We will apply \cref{thm:integ-noise-time-scaling}.
We have $ \rho(z) = z $, $ U(z) = \frac{1}{2} z^2 $, and $ \phi_1(\theta) = \cos(\theta) $.
Thus, we know $ V = |\nabla_z U|^2 - \nabla_z^2 U = z^2 - 1 $
so $ V + 2 = z^2 + 1 \geq 0 $
and $ |\rho|^2 = z^2 \leq z^2 + 1 $.
Let $ \{H_n\}_{n=0}^\infty $ be the Hermite polynomials.
$ A_{\Theta,Z}^{-1} e^{ik\theta} $ can be calculated by
expanding the sought function in terms of Hermite polynomials,
substituting into the equation,
and solving for the coefficients of the expansion.
This is standard so we just give the result.
Recall that $ H_n' = 2zH_n - H_{n+1} = 2n H_{n-1} $
so $ H_n'' - 2zH_n' = 2n (H_{n-1}' - 2z H_{n-1}) = -2n H_n $.
Now, we define $ \Phi_k : \mathbb{R} \times \mathbb{R} \to \mathbb{C} $ as
\[
\Phi_k(\theta, z) := e^{3k^2/4} e^{ik(\theta + z)} \sum_{n=0}^\infty \frac{(-ik)^n}{(n + k^2 / 2) 2^n n!} H_n(z - ik).
\]
The infinitesimal generator for $ (\Theta, Z) $ (the driving process without time scaling)
is $ A_{\Theta,Z} = \frac{1}{2} \partial_z^2 - z \partial_z + z \partial_\theta $ so
\[
e^{-ik(\theta + z)} A_{\Theta,Z} e^{ik(\theta + z)}
= \frac{1}{2} (\partial_z + ik)^2 - z (\partial_z + ik) + z(\partial_\theta + ik)
= \frac{1}{2} \partial_z^2 - (z - ik) \partial_z - \frac{k^2}{2} + z \partial_\theta
\]
and $ \left(\frac{1}{2} \partial_z^2 - (z - ik) \partial_z\right) H_n(z - ik) = -n H_n(z - ik) $.
Thus, using the generating function for the Hermite polynomials, namely $ e^{2zt - t^2} = \sum_{n=0}^\infty H_n(z) \frac{t^n}{n!} $,
$ A_{\Theta,Z} $ will act on $ \Phi_k(\theta, z) $ as
\[ \begin{split}
A_{\Theta,Z} \Phi_k
&= e^{3k^2/4} e^{ik(\theta + z)} \sum_{n=0}^\infty \frac{(-ik)^n}{(n + k^2 / 2) 2^n n!} \left(\frac{1}{2} \partial_z^2 - (z - ik) \partial_z - \frac{k^2}{2} + z \partial_\theta\right) H_n(z - ik) \\
&= - e^{3k^2/4} e^{ik(\theta + z)} \sum_{n=0}^\infty \frac{(-ik)^n}{2^n n!} H_n(z - ik)
= - e^{3k^2/4} e^{ik(\theta + z)} e^{2(z - ik)(-ik/2) - (-ik/2)^2}
= - e^{ik\theta}.
\end{split} \]
Therefore, $ - \Phi_k = A_{\Theta,Z}^{-1} e^{ik\theta} $.
We will let $ \mathcal{H} = L^2\left(\mathbb{R}, \frac{1}{\sqrt{\pi}} e^{-z^2} \dr z\right) $ which is the space in which $ \left\{(2^n n!)^{-1/2} H_n \right\} $ is an orthogonal basis.
Then, for any $ k, \ell \in \mathbb{Z} $,
we can calculate the value of the form
\[ \begin{split}
\langle e^{i \ell \theta}, e^{i k \theta} \rangle_{\Theta,Z}
&= - \int_{-\infty}^\infty \int_0^{2\pi}
e^{-i \ell \theta} (- \Phi_k(\theta, z))
\frac{1}{2\pi} \frac{e^{-z^2}}{\sqrt{\pi}} \dr \theta \dr z \\
&= e^{3k^2/4} \sum_{n=0}^\infty \frac{(-ik)^n}{(n + k^2/2) 2^n n!}
\left(\frac{1}{2\pi}\int_0^{2\pi} e^{i(k -\ell) \theta} \dr m\right)
\left(\frac{1}{\sqrt{\pi}} \int_{-\infty}^\infty e^{-ik(z - ik)} H_n(z - ik) e^{-z^2 + 2ikz + k^2} \dr z\right) \\
&= e^{3k^2/4} \delta_{k \ell} \sum_{n=0}^\infty \frac{(-ik)^n}{(n + k^2/2) 2^n n!} \frac{1}{\sqrt{\pi}} \int_{-\infty - ik}^{\infty - ik}  e^{-ikz} H_n(z) e^{-z^2} \dr z
\end{split} \]
By the holomorphy of $ e^{- z^2 - ikz} H_n(z) $
and because of the rapid decay of $ e^{-z^2} $, we can use the generating function and orthogonality of the Hermite polynomials to get
\[ \begin{split}
\frac{1}{\sqrt{\pi}} \int_{-\infty - ik}^{\infty - ik} e^{2 \left(-\frac{ik}{2}\right) z} H_n(z) e^{-z^2} \dr z
&= e^{\left(-\frac{ik}{2}\right)^2} \sum_{j=0}^\infty \frac{(-ik)^j}{2^j j!} \int_{-\infty}^\infty H_j(z) H_n(z) \frac{1}{\sqrt{\pi}} e^{-z^2} \dr z \\
&= e^{-k^2/4} \sum_{j=0}^\infty \frac{(-ik)^j}{2^j j!} 2^n n! \delta_{nj}
= e^{-k^2/4} (-ik)^n
\end{split} \]
and so $ \langle e^{i\ell \theta}, e^{ik \theta} \rangle_{\Theta,Z} $ becomes
\[
\langle e^{i \ell \theta}, e^{ik \theta} \rangle_{\Theta,Z}
= e^{3k^2/4} \delta_{k\ell} \sum_{n=0}^\infty
\frac{(-ik)^n}{(n + k^2/2) 2^n n!} e^{-k^2/4} (-ik)^n
= e^{\frac{k^2}{2}} \delta_{k\ell} \sum_{n=0}^\infty \frac{(-1)^n k^{2n}}{(n + k^2/2) 2^n n!}
\]
completing the proof.
\end{proof}

\textbf{Acknowledgements:}
J. W. thanks D. Bernard for the reference to \cite{BBJ}
which provided the insight, crucial for the present paper.  
Collaboration with G. Volpe led to an earlier analysis of light-sensitive robots, made rigorous here.
The content of \cref{ssub:mips} is a rigorous version of the results
obtained in a discussion with J. Stenhammar and G Volpe.
Discussions with J. Birrell, K. Burdzy, M. Latifi-Jebelli, J. Watkins, and E. Waymire are gratefully acknowledged.

\textbf{Funding, Declarations, and Conflicts of Interest:}
The research presented here was partially supported by a mini-grant from the University of Arizona Department of Mathematics, NSF grant DMS 1911358,
and by the Simons Foundation Fellowship 23539.
Besides these funding sources,
there are no further relevant financial or non-financial interests to disclose.
The authors have no conflicting interests
that are relevant to the content of this article.
There were no human or animal participants involved in the creation of this article
and there is no associated data.

\filbreak
\providecommand{\bysame}{\leavevmode\hbox to3em{\hrulefill}\thinspace}
\providecommand{\MR}{\relax\ifhmode\unskip\space\fi MR }
\providecommand{\MRhref}[2]{%
  \href{http://www.ams.org/mathscinet-getitem?mr=#1}{#2}
}
\providecommand{\href}[2]{#2}

\end{document}